\documentclass[a4paper]{amsart}

\usepackage{color}
\usepackage{graphicx,epstopdf}
\usepackage{amsmath}
\usepackage{amssymb}
\usepackage{amsfonts}
\usepackage{amsthm}

\theoremstyle{remark}
\newtheorem{remark}{Remark}[section]
 
\begin{document}
 \title[]{Multi-domain spectral approach to rational-order fractional 
 derivatives}

\author{Christian Klein$^{*}$}
\address{Institut de Math\'ematiques de Bourgogne, UMR 5584\\
Institut Universitaire de France,\\
                Universit\'e de Bourgogne-Franche-Comt\'e, 9 avenue Alain Savary, 21078 Dijon
                Cedex, France\\
    E-mail Christian.Klein@u-bourgogne.fr}

\author{Nikola Stoilov}
\address{Institut de Math\'ematiques de Bourgogne, UMR 5584\\
               Universit\'e de Bourgogne-Franche-Comt\'e, 9 avenue Alain Savary, 21078 Dijon
               Cedex, France\\
   E-mail Nikola.Stoilov@u-bourgogne.fr}
\date{\today}

\begin{abstract}
We propose a method to numerically compute fractional derivatives (or 
the fractional Laplacian) on the whole real line via Riesz fractional integrals. The compactified real line is 
divided into a number of intervals, thus amounting to a multi-domain 
approach; after transformations in accordance with the underlying 
$Z_{q}$ curve ensuring analyticity of the respective integrands, 
the integrals over the different domains are computed with a Clenshaw-Curtis algorithm. As an example, we consider solitary waves for fractional 
Korteweg-de Vries equations and  compare these to results 
obtained with a discrete Fourier transform. 
 
\end{abstract}

 
\thanks{This work was partially supported by 
 the ANR-17-EURE-0002 EIPHI and by the 
European Union Horizon 2020 research and innovation program under the 
Marie Sklodowska-Curie RISE 2017 grant agreement no. 778010 
IPaDEGAN}
\maketitle

\section{Introduction}

Fractional derivatives have gained considerable interest in 
applications. They  appear for instance when dealing with 
flows in heterogeneous media \cite{BWM00}, modeling viscoelasticity 
in biological and medical applications \cite{AEL08, YPK16}, or 
anomalous diffusion \cite{KS05, YAL04} and financial flows 
\cite{CT04}. An important aspect is that in contrast to 
classical derivatives, they are
not local.  Temporal fractional derivatives thus allow to include memory 
effects in partial differential 
equations (PDEs), and fractional spatial derivatives allow to model 
nonlocal effects. In the context of  dispersive PDEs, nonlocal effects are particularly interesting in the theory of 
water waves, see \cite{Lan} for a comprehensive review. The motion of 
a surface of water, say the sea, under the influence of gravity leads 
to a nonlocal equation. The \emph{Whitham equation} is a nonlocal 
equation in 1D with the exact dispersion of the water wave equations, 
see \cite{KLPS} for references. An interesting toy model in this 
context  is  the fractional Korteweg-de 
Vries (KdV) equation,
\begin{equation}
	u_{t}-D^{\alpha}u_{x}+uu_{x}=0
	\label{fKdV},
\end{equation}
which allows to study the interaction of dispersion and nonlinearity in PDEs, see 
\cite{KSbook} for a recent review with many references;
here the fractional derivative $D^{\alpha}$ can be defined 
via its Fourier symbol
\begin{equation}
	\mathcal{F}D^{\alpha} = |k|^{\alpha},\quad 
	0<\alpha<1,
	\label{fourier}
\end{equation}
where $\mathcal{F}$ denotes the Fourier transform in $x$, and where 
$k$ is the dual Fourier variable, see (\ref{fourierdef}) for a 
definition. These fractional 
derivatives are sometimes referred to as the fractional 
Laplacian $(-\Delta)^{\alpha/2}$ (for different definitions of 
fractional derivatives see for instance \cite{Hil} and references 
therein). 

In this paper we are interested in the efficient computation of fractional 
derivatives on the whole real line. The fractional derivative 
(\ref{fourier}) can be also computed in terms 
of Riesz fractional integrals over the whole real line, see for 
instance \cite{Hil}, 
\begin{equation}
	D^{\alpha}u(x) = 
	\frac{1}{2\Gamma(1-\alpha)\sin((1-\alpha)\pi/2)}\partial_{x}\int_{-\infty}^{\infty} 
	\frac{\mbox{sgn}(x-y)u(y)}{|x-y|^{\alpha}}dy
	\label{fracint}.
\end{equation}
In this work we present a multi-domain spectral approach to 
efficiently compute these integrals for general values of $x\in 
\mathbb{R}\cup\{\infty\}$. Multi-domain methods have the advantage to 
combine the main feature of spectral methods, the approximation of a 
smooth function with a numerical error decreasing exponentially with 
the resolution, with a concept of locality (note that spectral methods 
are always global on the considered interval). They are similar to 
\emph{spectral elements} in combining high order convergence with the 
freedom to choose the intervals. Thus resolution can be attributed 
where needed, for instance for highly oscillatory zones of the 
studied function. An important point is that 
different forms of the same equation or integral can be studied in different 
domains which allows to use adapted variables, for instance in 
the context of a Riemann surface underlying fractional derivatives,
see \cite{RSbook} and references. These methods have been successfully 
applied to various  problems in the context of differential equations 
as in \cite{CFKSV} for the hypergeometric equation, in \cite{KRS} for 
the Hilbert transform and in
\cite{KPS} for Fourier transforms in the context of dispersive PDEs.

 
In the present paper, the real line is subdivided into a 
finite number of intervals two of which are unbounded. On each of the 
considered intervals, coordinate transformations are applied to 
obtain a smooth integrand that is numerically computed with the 
Clenshaw-Curtis algorithm \cite{CC}. This algorithm shows as 
mentioned \emph{spectral convergence}, the numerical error in integrating a smooth 
function will decrease exponentially with the numerical resolution.

The algorithm is tested for known examples and then applied 
to the solitary waves of the fractional KdV equation (\ref{fKdV}), 
i.e., solutions of the form $u(x,t)=Q_{c}(x-ct)$ vanishing at infinity where $c>0$ is a real 
constant. These solitary waves satisfy the equation 
\begin{equation}
	cQ_{c}+D^{\alpha}Q_{c}-\frac{1}{2}Q_{c}^{2}=0
	\label{fKdVsol}.
\end{equation}
We numerically construct solutions of this equation with a 
Newton-Krylov method and and compare them to 
solutions previously constructed with a discrete Fourier transform 
(DFT) in \cite{KS15}.

The paper is organized as follows: in section 2 we collect some known 
facts on fractional derivatives and algebraic curves. In section 3 we 
present the multi-domain approach to compute the 
integral (\ref{fracint}). In section 4 we outline the numerical 
approach and test it for known examples. The results are compared to 
what can be achieved with DFT techniques. In section 5 we construct 
solitary waves for the fractional KdV equation. We add some 
concluding remarks in section 6. 

\section{Overview and Preliminaries}

In this section we present a brief overview of existing numerical approaches to fractional derivatives and collect some facts on the subject that will be   
needed in the following. 

\subsection{Numerical approaches}

A number of different approaches to numerically computing fractional derivatives are actively worked on. 
These include adaptive finite element methods that use much finer 
grids close to the boundary of integration domains, to address the 
lower regularity of fractional derivatives close to said boundaries 
\cite{AG17} and finite difference methods on stencils in the complex 
plane \cite{FP}. 
Walking on spheres is a classical Monte Carlo method, recently 
applied to computing the fractional Poisson equation \cite{KOT18}, memory-efficient approaches based on convolution quadrature via Runge-Kutta schemes \cite{BLF19}; an excellent review can be found in  \cite{LetAL20}

We should mention using various orthogonal bases, for example 
ultraspherical polynomials \cite{HO18}, Jacobi 
polynomials in fractional powers \cite{PF}, as well as methods based on Laguerre polynomials \cite{CSW18}, and using 
specially selected super-convergence points \cite{DZZ19}. 
Interesting approaches to fractional integrals are based on a 
result by Bateman \cite{Bat} that certain hypergeometric functions are 
eigenfunctions of the fractional derivatives. This was used in 
\cite{CCH} to compute fractional integrals in terms of 
hypergeometric functions. Note that we used in \cite{CFKSV} a similar 
approach as in the present paper to compute hypergeometric functions 
in the whole complex plane. 

\subsection{Fractional integrals}

To fix notation, we use the standard definition for a  Fourier transform for 
functions $u(x)$ 
in the Schwartz class of rapidly decreasing smooth functions 
denoted by $(\mathcal{F} u)(k)$ with dual variable $k$ 
and its inverse,
\begin{align}
	(\mathcal{F}u)(k)& = \int_{\mathbb{R}}^{}u(x)e^{-ikx}dx,\quad k\in 
	\mathbb{R},
	\nonumber\\
	u(x) & =\frac{1}{2\pi}\int_{\mathbb{R}}^{}(\mathcal{F}u)(k)e^{ikx} 
	dk,\quad x\in\mathbb{R}.
	\label{fourierdef}
\end{align}
The convention extends by duality to tempered distributions.

The 
conjugate Riesz fractional  integral of order $\alpha$ reads (see for 
instance \cite{Hil})
\begin{equation}
	\widetilde{I^{\alpha}u}(x)=\frac{1}{2\Gamma(\alpha)\sin(\alpha\pi/2)}\int_{-\infty}^{\infty} 
	\frac{\mbox{sgn}(x-y)u(y)}{|x-y|^{1-\alpha}}dy
	\label{Itilde}.
\end{equation}
For the Fourier transform of this integral, one has
\begin{equation}
	\mathcal{F}(\widetilde{I^{\alpha}u}(x)) = 
	-i\mbox{sgn}(k)|k|^{-\alpha}(\mathcal{F}u)(k).
	\label{Itildef}
\end{equation}
Thus we get 
\begin{equation}
	\mathcal{F}(\partial_{x}\widetilde{I^{1-\alpha}u}(x))(k) = 
	|k|^{\alpha}(\mathcal{F}u)(k).
	\label{fracdef}
\end{equation}
This expression leads to (\ref{fracint}) for the fractional 
derivative to be studied in this paper. 
Note that there is an equivalent definition of these fractional 
derivatives (or fractional Laplacians $(-\Delta^{\alpha/2})$, see for instance \cite{CCH}),
\begin{equation}
	D^{\alpha}u = \frac{\alpha 2^{\alpha-1} \Gamma(1/2+\alpha/2)}{ 
		\sqrt{\pi}\Gamma(1-\alpha/2)}	\int_{\mathbb{R}}^{}\frac{u(x)-u(x+y)}{|y|^{1+\alpha}}dy
	\label{Lap}.
\end{equation}
We do not use the latter version of the fractional integrals here 
since the integrand would be singular without the inclusion of the 
term $u(x)$ whereas this is not the case in (\ref{fracint}). However 
the form (\ref{Lap}) can be treated with a combination of the methods 
detailed in this paper and in \cite{KRS}. The price to pay to apply 
the integral (\ref{fracint}) is the appearance of a classical 
derivative. However since we are interested in using these 
derivatives for iterative approaches and thus recurrent computations 
of the fractional derivatives, we anyway work on numerical grids 
where it is straightforward to compute classical derivatives. 

This means we will have to compute integrals of the form
\begin{equation}
	I^{-}(x)=\int_{-\infty}^{x}\frac{u(y)}{(x-y)^{\alpha}}dy,\quad 
	I^{+}(x)=\int_{x}^{\infty}\frac{u(y)}{(y-x)^{\alpha}}dy
	\label{Ipm},
\end{equation}
where $0<\alpha<1$, and where $u:\mathbb{R}\mapsto \mathbb{R}$. We 
assume that $\alpha$ is a rational number, and we are interested 
in this integral for $x\in \mathbb{R}$. 

\subsection{$Z_{q}$ curves}

Note that the integrals (\ref{Ipm}) are Abelian integrals defined on a Riemann 
surface. If $\alpha$ is not rational, this Riemann surface will not be 
compact. Here we consider always the case of rational $\alpha$  with 
$0<\alpha<1$, i.e., 
\begin{equation}
	\alpha = \frac{p}{q},\quad p,q\in \mathbb{N},\quad p<q,
	\label{pq}
\end{equation}
where $p$, $q$ are coprime. This ensures that the underlying Riemann 
surfaces are compact. Thus the integrals (\ref{Ipm}) 
are integrals on the surface defined 
by the family of $Z_{q}$ curves (parametrised by the branch point $x$)
\begin{equation}
	\mu^{q}=K-x, \quad \mu,K\in \mathbb{C}
	\label{Zq}
\end{equation}
branched at $x$
and infinity, see for instance \cite{RSbook} for references and 
computational approaches to such surfaces (general $Z_{q}$ curves are 
given by the algebraic equation $\mu^{q}$ equals a polynomial in 
$K$). This means the fractional integrals (\ref{Ipm}) are functions 
on the moduli space of the $Z_{q}$ curves (\ref{Zq}). Analytic 
functions on these curves will not be analytic functions in $K$ near 
the branch points, however, they are analytic functions in the local parameters on 
the curve: we have $\lambda_{x}=(K-x)^{1/q}$ near $x$ and 
$\lambda_{\infty}=1/K^{1/q}$ near infinity. Series expansions in terms of these 
local parameters and thus fractional powers of $K$ are called 
\emph{Puiseux series}, see again  
\cite{RSbook}. In this paper 
near infinity we use the parameters $\xi_{\pm}:=(\pm1/x)^{1/q}$ which allows us 
to also consider functions having a cusp for $1/x=0$.  Obviously there is only 
one local parameter needed in the vicinity of infinity, but since we 
are working on $\mathbb{R}$ only, we introduce different ones near 
$\pm \infty$.

 It is the goal of this paper to actively 
exploit the underlying Riemann surface structure in the computation 
of the fractional integrals (\ref{Ipm}). 
Note that the 
differentials $dK/\mu^{p}$ in (\ref{Ipm}) are not necessarily holomorphic, i.e., 
of the form $f(z)dz$ with $f(z)$ holomorphic in any chart on the underlying Riemann surface. 
A fall off condition at infinity is needed for the 
integrals to exist for general $q>p>0$. We assume in the following the fall off condition 
of  solitary waves of the fKdV and fractional nonlinear Schr\"odinger 
equation as proven by Frank 
and Lenzmann \cite{FL}, 
\begin{equation}
	\lim_{|x|\to \infty}u(x)|x|^{1+\alpha}<\infty
	\label{falloff}.
\end{equation}
Note that this condition can be slightly relaxed, all we need for our 
numerical approach to show spectral convergence (an exponential 
decrease of the numerical error with the resolution for analytic 
functions) is that the differential is holomorphic in a neighborhood near 
infinity. This is the case if $\lim_{|y|\to\infty}u(y)|y|$ is finite. 
The code can be applied to this case, but we only discuss here the 
fall off condition of fractional KdV solitons. 
Thus we define 
\begin{align}
		u^{I}(\xi_{-}) &:= u(x)(-x)^{1+\alpha}, \quad x<a, \\
		u^{II}(x) &:= u(x), \quad  a\leq x\leq b, \\
		u^{III}(\xi_{+})& := u(x)x^{1+\alpha},\quad  x>b.
	\label{cases}
\end{align}
\begin{remark}
	If the function $u$ is not analytic on $\mathbb{R}$, the 
	resulting fractional derivative will be singular at the 
	discontinuities. This can be 
	seen for an example with compact support, say 
$u$ being the characteristic function of the interval $[a,b]$. One 
gets for (\ref{fracint})
\begin{equation}
	D^{\alpha}u(x) = 
	\frac{1}{2\Gamma(1-\alpha)\sin((1-\alpha)\pi/2)}
\left(\frac{1}{(x-a)^{\alpha}}+\frac{1}{(b-x)^{\alpha}}\right)
	\label{char}.
\end{equation}
 If $\lambda_{a}$, 
$\lambda_{b}$ are the local parameters on the $Z_{q}$ curve near the 
branch points $a$, $b$, one has $D^{\alpha}u(x) = 
\frac{1}{2\Gamma(1-\alpha)\sin((1-\alpha)\pi/2)}(\lambda_{a}^{-p}+\lambda_{b}^{-p})$. 
\end{remark}
Thus in a multi-domain approach as in the next section, one can simply divide the interval 
into two intervals, say $[a,(a+b)/2]$ and $[(a+b)/2,b]$, and compute 
the integrals $I^{\pm}$ using the local parameter near a branch point 
on the curve. Then the classical derivative with respect to 
$\partial_{\lambda_{a}}$ is computed near $a$
giving the wanted result up to a factor $\lambda_{a}^{-p}$ and 
similarly near $b$. The singular factor can be taken care of 
explicitly. To avoid related problems, we will only consider fractional 
derivatives of smooth functions in this paper, but it is straightforward to generalize the approach to piecewise smooth functions.

\section{Multi-domain approach}
In this section we discuss how to set up a multi-domain approach
for the computation of the fractional 
derivative via the integral (\ref{fracint}). The basic idea is to use 
the local parameters of the underlying $Z_{q}$ curve in the integrals. 
To this end we
divide the real line into at least three domains: we choose two real 
constants $a,b$ with $a<0<b$. Then domain I  is given by $x<a<0$, domain II 
by $a\leq x \leq b$, and domain III by $x>b>0$. Note that there can be 
any number of finite domains (type II). We discuss only the case of 
one such domain here in order to avoid cluttered notation, but the 
generalization to more finite domains is straightforward.

We write the fractional integral in (\ref{fracint}) as 
$$\frac{1}{2\Gamma(1-\alpha)\sin((1-\alpha)\pi/2)}\partial_{x} 
(I^{-}(x)-I^{+}(x)),$$ where $I^{\pm}$ are defined in (\ref{Ipm}). 
The integrals $I^{-}$ and $I^{+}$ are split into several integrals 
for which substitutions depending on the value of $x$ are chosen in 
order that the integrands are smooth functions in all cases. We 
detail this for $I^{-}$ below, the ones for $I^{+}$ follow 
analogously.

\paragraph{\textbf{Domain I:}}
For the first interval $x<a<0$, we use as mentioned the local 
parameter  
$\xi_{-}=1/(-x)^{1/q}$ and also the substitution $s = -1/y$. Thus we get for $I^{-}$ in (\ref{Ipm})
\begin{equation}
		I^{-}(\xi_{-})=\xi_{-}^{p}\int_{0}^{\xi_{-}^{q}}\frac{u^{I}(s^{1/q})}{(\xi_{-}^{q}-s)^{\alpha}}s^{2\alpha-1}ds
	= I^{-}_{1}+I^{-}_{2}
	\label{I}.
\end{equation}
The integrand has possible singularities for $s=0$ and 
$s=\xi_{-}^{q}$. To address these, we separate the integral into two, 
$I_{1}^{-}$ from $0$ to $\xi_{-}^{q}/2$ and $I^{-}_{2}$ from  $\xi_{-}^{q}/2$ to  
$\xi_{-}^{q}$. For the first integral we apply the mapping 
$s=t^{q}$ (the local parameter near infinity of the $Z_{q}$ curve) to obtain
\begin{equation}
	I^{-}_{1}=	
	\xi_{-}^{p}\int_{0}^{\xi_{-}^{q}/2}\frac{u^{I}(s^{1/q})}{(\xi_{-}^{q}-s)^{\alpha}}s^{2\alpha-1}ds=
	q\xi_{-}^{p}\int_{0}^{\xi_{-}/2^{1/q}}\frac{u^{I}(t)}{(\xi_{-}^{q}-t^{q})^{\alpha}}t^{2p-1}dt\\
	\label{I1}.
\end{equation}
Note that $p\in \mathbb{N}$ which means that the integrand is smooth 
for $t\to0$. 
For the second integral, we apply $s=\xi_{-}^{q}-t^{q}$ (the local 
parameter near $x$) to get once 
more a smooth integrand, 
\begin{equation}
	I^{-}_{2}=	
	\xi_{-}^{p}\int_{\xi_{-}^{q}/2}^{\xi_{-}^{q}}\frac{u^{I}(s^{1/q})}{(\xi_{-}^{q}-s)^{\alpha}}s^{2\alpha-1}ds=
	q\xi_{-}^{p}\int_{0}^{\xi_{-}/2^{1/q}}u^{I}((\xi_{-}^{q}-t^{q})^{1/q})(\xi_{-}^{q}-t^{q})^{2\alpha-1}t^{q-p-1}dt.
	\label{I2}
\end{equation}

\paragraph{\textbf{Domain II:}}
For $x\in [a,b]$, we split the integral $I^{-}$ once more into two 
integrals, $I^{-}(x) = I^{-}_{3}+I^{-}_{4}$ where 
\begin{equation}
	I^{-}_{3} = 	\int_{a}^{x}\frac{u(y)}{(x-y)^{\alpha}}dy=
	q\int_{0}^{(x-a)^{1/q}}u^{II}(x-t^{q})
	t^{q-p-1}dt,
	\label{I3}
\end{equation}
where we have put $y=x-t^{q}$, and where 
\begin{equation}
	I^{-}_{4}=\int_{-\infty}^{a}\frac{u(y)}{(x-y)^{\alpha}}dy=q\int_{0}^{(-1/a)^{1/q}}\frac{u^{I}(
	t)}{(1+xt^{q})^{\alpha}}t^{2p-1}dt
	\label{I4},
\end{equation}
where we have put $y = -1/t^{q}$ as before.
If $x$ is close to $a$, i.e., $x<a+\delta$ with $0<\delta\ll1$, 
the  integral in (\ref{I4}) will have an almost singular 
integrand which will affect the convergence of a spectral method. 
In this case we will apply the same approach as in 
(\ref{I}), i.e.,
\begin{align}
	I^{-}_{4}  &= 
		q\int_{0}^{(-1/(2a))^{1/q}}\frac{u^{I}(t)}{(1+xt^{q})^{\alpha}}t^{2p-1}dt
\nonumber\\
	& +\frac{q}{x}\int_{(1-x/(2a))^{1/q}}^{(1-x/a)^{1/q}}
	u^{I}(((t^{q}-1)/x)^{1/q})\left(\frac{t^{q}-1}{x}\right)^{2\alpha-1} 
	t^{q-p-1}dt
	\label{lint22}.
\end{align}

\paragraph{\textbf{Domain III:}}
In the third interval $x>b$, the fractional derivative is computed as 
$
	I^{-} = I_{5}^{-}+I^{-}_{6}+I^{-}_{7}
$ with a local variable $\xi_+ = 1/x^{1/q}$
where 
\begin{equation}
	I^{-}_{5}=
	\int_{-\infty}^{a}\frac{u(y)}{(x-y)^{\alpha}}dy = q\xi_{+}^{p}
	\int_{0}^{(-1/a)^{1/q}}\frac{u^{I}(t)}{(\xi_{+}^{q}+t^{q})^{\alpha}}
	t^{2p-1}dt,
	\label{largea}
\end{equation}
 where 
\begin{equation}
	I^{-}_{6} = \int_{a}^{b}\frac{u(y)}{(x-y)^{\alpha}}dy = 
	\xi_{+}^{p}\int_{a}^{b}\frac{u^{II}(y)}{(1-\xi_{+}^{q}y)^{\alpha}}dy
	\label{I6},
\end{equation}
and where 
\begin{equation}
	I^{-}_{7} = \int_{b}^{x}\frac{u(y)}{(x-y)^{\alpha}}dy = 
	q\xi_{+}^{p}\int_{0}^{(1/b-\xi_{+}^{q})^{1/q}}
	u^{III}((\xi_{+}^{q}+t^{q})^{1/q})(\xi_{+}^{q}+t^{q})^{2\alpha-1} t^{q-p-1}dt
	\label{large2}.
\end{equation}
If $\xi_{+}<1/b^{1/q}-\delta$, we compute $I^{-}_{6}$ via
\begin{equation}
	I_{6}^{-} = 
	\int_{(1-\xi_{+}^{q}a)^{1/q}}^{(1-\xi_{+}^{q}b)^{1/q}}u^{II}\left(\frac{1-t^{q}}{\xi_{+}^{q}}\right) t^{q-p-1}dt.
	\label{Im6}
\end{equation}
If $\xi_{+}$ is small, one is close to the branch 
point of the root in (\ref{large2}) in the argument of $u^{III}$. In 
such a case, the root is not well approximated by polynomials. Therefore 
we split for $\xi_{+}<\delta$ the integral $I_{7}^{-}$ into two 
integrals and use the above substitutions to get regular integrands, 
\begin{equation}
	\begin{split}
	I_{7}^{-}&= 
	q\xi_{+}^{p}\int_{0}^{(\gamma^{q}-1)^{1/q}\xi_{+}}u^{III}((t^{q}+\xi_{+}^{q})^{1/q})(\xi_{+}^{q}+t^{q})^{2\alpha-1} t^{q-p-1}dt 
	\\
	&+q\xi_{+}^{p}\int_{\gamma\xi_{+}}^{1/b^{1/q}}\frac{u^{III}(t)}{(t^{q}-\xi_{+}^{q})^{\alpha}}
	t^{2p-1}dt,
	\end{split}
	\label{I70}
\end{equation}
where $\gamma$ is some constant greater than 1 (we take the minimum of 
$1/\xi_{+}/2$ and $1/\delta$). For the same reason we write 
$I_{5}^{-}$ (\ref{largea}) in the form (essentially just splitting 
the integral into two integrals)
\begin{equation}
	I^{-}_{5}= q\xi_{+}^{p}
	\int_{0}^{\gamma\xi_{+}}\frac{u^{I}(t)}{(\xi_{+}^{q}+t^{q})^{\alpha}}
	t^{2p-1}dt+\ q\xi_{+}^{p}
	\int_{\gamma\xi_{+}}^{(-1/a)^{1/q}}\frac{u^{I}(t)}{(\xi_{+}^{q}+t^{q})^{\alpha}}
	t^{2p-1}dt.
	\label{I50}
\end{equation}

\begin{remark}
	The integrals $I_{\pm}$ are in principle so-called periods of 
	some differentials on the 
	Riemann surface defined by the $Z_{q}$-curve (\ref{Zq}). This 
	means integrals along cycles on the Riemann surface around the 
	branch points $x$ and infinity and thus 
	non-homologous to zero. These periods can be reduced to integrals 
	along the cut between the branch points. Numerically it would be 
	convenient to avoid the branch points by integrating on contours 
	in the complex plane, see for instance \cite{RSbook}, but to this 
	end one would need the function $u$ in the complex plane and its 
	singularities. Since we consider $u$  here 
	only on the real axis, one has 
	to impose the above substitutions. 
\end{remark}

\paragraph{\textbf{The integral $I^{+}$}}

In a similar way we want to compute the integral $I^{+}$ in 
(\ref{Ipm}). Essentially one has to interchange domains I and III to 
obtain the wanted form of the integrals from the above relations for 
$I^{-}$. In 
domain III ($x>b$), we write $I^{+}=I^{+}_{1}+I^{+}_{2}$ with  
\begin{equation}
	I^{+}_{1}=	
	q\xi_{+}^{p}\int_{0}^{\xi_{+}/2^{1/q}}\frac{u^{III}(t)}{
	(\xi_{+}^{q}-t^{q})^{\alpha}}t^{2p-1}dt
	\label{Ip1}
\end{equation}
and 
\begin{equation}
	I^{+}_{2}=	q\xi_{+}^{p}\int_{0}^{\xi_{+}/2^{1/q}}
	u^{III}((\xi_{+}^{q}-t^{q})^{1/q})(\xi_{+}^{q}-t^{q})^{2\alpha-1}t^{q-p-1}dt.			
	\label{Ip2}
\end{equation}

In domain II, $a<x<b$ we write $I^{+}=I_{3}^{+}+I_{4}^{+}$ with 
\begin{equation}
	I^{+}_{3}=
	q\int_{0}^{(1/b)^{1/q}}\frac{u^{III}(t)}{
	(1-xt^{q})^{\alpha}}t^{2p-1}dt	
	\label{IIp3}
\end{equation}
and with 
\begin{equation}
	I^{+}_{4}=q\int_{0}^{(b-x)^{1/q}}
	u^{II}((x+t^{q}))t^{q-p-1}dt
	\label{IIp4}.
\end{equation}
For $x$ near $b$, $x>b-\delta$, we use instead of (\ref{IIp3})
\begin{equation}
	\begin{split}
	I^{+}_{3}&=		
	q\int_{0}^{(1/(2b))^{1/q}}\frac{u^{III}(t)}{(1-xt^{q})^{\alpha}}t^{2p-1}dt\nonumber \\
	&-\frac{q}{x}\int_{(1-x/(2b))^{1/q}}^{(1-x/b)^{1/q}}u^{III}(((1-t^{q})/x)^{1/q})\left(\frac{1-t^{q}}{x}\right)^{2\alpha-1} 
	t^{q-p-1}dt.
	\end{split}
	\label{IIp32}
\end{equation}

In domain I, $x<a$, we write $I^{+}=I^{+}_{5}+I^{+}_{6}+I^{+}_{7}$ 
with 
\begin{equation}
	I^{+}_{5}=q\xi_{-}^{p}\int_{0}^{(1/b)^{1/q}}\frac{u^{III}(t)}{
	(\xi_{-}^{q}+t^{q})^{\alpha}}t^{2p-1}dt
	\label{Ip5},
\end{equation}
with 
\begin{equation}
	I^{+}_{6}=\xi_{-}^{p}\int_{a}^{b}\frac{u^{II}(y)}{(\xi_{-}^{q}y+1)^{\alpha}}dy 
	\label{Ip6}
\end{equation}
and with 
\begin{equation}
	I^{+}_{7}=q\xi_{-}^{p}\int_{0}^{(-\xi_{-}^{q}-1/a)^{1/q}}u^{I}((\xi_{-}^{q}+t^{q})^{1/q})
	(\xi_{-}^{q}+t^{q})^{2\alpha-1}t^{q-p-1}dt	
	\label{Ip7}.
\end{equation}
For $\xi_{-}>1/(-a)^{1/q}-\delta$, the integral (\ref{Ip6}) is computed in the form
\begin{equation}
	I^{+}_{6}=q(\xi_{-})^{p-q}
	\int_{(1+a\xi_{-}^{q})^{1/q}}^{(1+b\xi_{-}^{q})^{1/q}}u^{II}\left(\frac{t^{q}-1}{\xi_{-}^{q}}\right) 
	t^{q-p-1}dt.
	\label{I2a}
\end{equation}
For $\xi_{-}<\delta$, we use
\begin{equation}
	\begin{split}
	I_{7}^{+}&= 
	q\xi_{-}^{p}\int_{0}^{(\gamma^{q}-1)^{1/q}\xi_{-}}u^{I}((t^{q}+\xi_{-}^{q})^{1/q})(\xi_{-}^{q}+t^{q})^{2\alpha-1} t^{q-p-1}dt \\
	&+q\xi_{-}^{p}\int_{\gamma\xi_{-}}^{1/(-a)^{1/q}}\frac{u^{I}(t)}{(t^{q}-\xi_{-}^{q})^{\alpha}}
	t^{2p-1}dt
	\end{split}
	\label{I70p},
\end{equation}
where $\gamma$ is some constant greater than 1 (we take the minimum of 
$1/\xi_{-}/2$ and $1/\delta$); in addition we write 
$I_{5}^{-}$ (\ref{largea}) in the form 
\begin{equation}
	I^{+}_{5}= q\xi_{+}^{p}
	\int_{0}^{\gamma\xi_{+}}\frac{u^{III}(t)}{(\xi_{+}^{q}+t^{q})^{\alpha}}
	t^{2p-1}dt+\ q\xi_{+}^{p}
	\int_{\gamma\xi_{+}}^{(1/b)^{1/q}}\frac{u^{III}(t)}{(\xi_{+}^{q}+t^{q})^{\alpha}}
	t^{2p-1}dt.
	\label{I50p}
\end{equation}

Note that the choices of the parameters $\delta$ and $\gamma$ in the 
relations above need not be the same for each case. We just want to 
avoid a cluttering of notation, and for the examples to be discussed 
in this paper, identical parameters can be applied to all cases.

\section{Numerical approach}
In the previous section, we have written the integrals $I^{\pm}(x)$ (\ref{Ipm}) 
for a given value of $x$ in terms of a sum of integrals. Each of these 
integrals has a smooth integrand and can be computed conveniently 
with the Clenshaw-Curtis algorithm \cite{CC}.  In this section, we 
explain how this will be done, and discuss some examples. The 
multi-domain approach is also compared to results obtained with a 
DFT. 

\subsection{Chebyshev polynomials}
The integrals of the previous section are all of the form 
\begin{equation}
	\int_{c}^{d}u(f(t)) g(t) dt
	\label{intform},
\end{equation}
where $c,d$ are some real constants, and where $f,g$ are some known 
smooth
functions. 

To compute the  integrals  (\ref{intform}), we map 
them to the interval $[-1,1]$ 
via  $t=d(1+l)/2+c(1-l)/2$, where $l\in[-1,1]$. 
The integrals of the form 
\begin{equation}
    \int_{-1}^{1}h(l)dl
    \label{gs}
\end{equation}
are computed with the Clenshaw-Curtis algorithm \cite{CC}: one
introduces the \emph{ Chebyshev collocation points }
\begin{equation}
    l_{m} = \cos(m\pi/N),\quad m = 0,\ldots,N
    \label{lm}.
\end{equation}
With some weights $w_{m}$, see \cite{trefethen} for a 
discussion and a code to compute them, the 
integral (\ref{gs}) is approximated via
\begin{equation}
     \int_{-1}^{1}h(l)dl\approx \sum_{m=0}^{N}w_{m}h(l_{m})
    \label{ccalg}
\end{equation}
which means that  just a 
scalar product has to be computed. The Clenshaw-Curtis algorithm 
is a spectral method,  see \cite{CC}
for an error analysis. 

In general the function $u$ will be only known on Chebyshev 
collocation points in each domain, for instance if one wants to 
determine solitary waves in an iteration. To obtain the function 
values on the collocation points in the variable $l$ on the interval 
$[-1,1]$, we use efficient and numerically stable \emph{barycentric 
interpolation}, see \cite{BT} for a review with many references. For convenience we compute the fractional integrals 
on Chebyshev collocation points on the 3 domains, but barycentric 
interpolation allows to obtain them for arbitrary values of $x\in 
\mathbb{R}$. The use of these Chebyshev collocation points also 
allows to compute the derivative in (\ref{fracint}) with 
\emph{Chebyshev differentiation matrices}, see \cite{trefethen,WR},
\begin{equation}
    \partial_{x} \vec{I} \approx D \vec{I}
    \label{gprime},
\end{equation}
where $\vec{I}$ is the vector with the components 
$I(l_{0}),\ldots,I(l_{N})$, on the Chebyshev 
collocation points.  Note that these differentiation matrices have a 
conditioning of the order of $\mathcal{O}(N^{2})$, see 
\cite{trefethen}. Thus it is generally recommended to avoid large 
values of $N$ since these will lead to a loss of accuracy. But in a 
multi-domain approach, this can be easily done by simply introducing 
more domains if needed. Note that these derivatives 
can be computed also with an FFT, see \cite{TW,trefethen}.

An advantage of the use of Chebyshev collocation points is that these 
can be related to an expansion of the considered function $f(x)$ in terms 
of Chebyshev polynomials \( 
T_{n}(x):=\cos(n\arccos(x)) \), 
\begin{equation}
	f(x)
\approx \sum_{n=0}^{N}c_{n}T_{n}(x)
	\label{chebcoeff}.
\end{equation}
These coefficients $c_{n}$, $n=0,1,\ldots,N$ can be 
computed efficiently via a \emph{fast cosine transform} which is 
related to the FFT, see e.g.~\cite{trefethen}. As for the DFT, the 
sum (\ref{chebcoeff}) can be seen as a truncated series, and the 
numerical error in doing so is indicated by the highest Chebyshev 
coefficients. Thus we can use the coefficients $c_{n}$ as a 
control of the numerical resolution in cases where no a priori estimates 
exist. For smooth functions, the coefficients decrease exponentially 
with the index as in the case of DFT coefficients, see the discussion 
in \cite{trefethen}. 

For the ease of presentation, we will use the same number $N$ of 
collocation points for all integrals though one could apply different 
resolutions. By checking the  coefficients we assure that the 
Chebyshev coefficients decrease to machine precision ($10^{-16}$ in 
our case, in practice limited to values of the order of $10^{-12}$ 
because of rounding errors) in all cases. For functions smooth on the whole real line, 
we adjust the interval limits $a$ and $b$ in a way that a similar 
resolution is needed for all integrals. Note that the computation for 
each value of $x$ is independent and thus fully parallelizable.

\subsection{Clenshaw-Curtis integration in the multi-domain approach}
In the previous section, the integrals $I^{\pm}$ in (\ref{fracint}) 
have been written in the various domains in terms of integrals with 
smooth integrands by taking advantage of the underlying Riemann surface. This 
means we introduced in the vicinity of the branch points of the 
$Z_{q}$ curve standard local coordinates in order to have 
smooth integrands. The latter will then be integrated with the 
Clenshaw-Curtis algorithm explained above. We illustrate this for the 
integral (\ref{I}) written in terms of (\ref{I1}) and (\ref{I2}). For 
the former, we consider the mapping $t=\xi_{-}/2^{1/q}(1+l)/2$ with 
$l\in[-1,1]$. Thus we get 
$$I^{-}_{1}=q\xi_{-}^{p}\xi_{-}/2^{1/q}/2\int_{-1}^{1}\frac{u^{I}(\xi_{-}/2^{1/q}(1+l)/2)}{
(\xi_{-}^{q}-(\xi_{-}/2^{1/q}(1+l)/2)^{q})^{\alpha}}dl.$$
This integral is of the form (\ref{gs}) with a smooth integrand and 
is approximated via (\ref{ccalg}). In a similar way we treat the 
integral (\ref{I2}) (same mapping for $t$) which leads to
\begin{eqnarray*}
I^{-}_{2}=q\xi_{-}^{p}\xi_{-}/2^{1/q}/2\int_{-1}^{1}u^{I}[\xi_{-}^{q}(1-1/2^{q}(1+l)^{q}/2)]
(\xi_{-}^{q}-(\xi_{-}/2^{1/q}(1+l)/2)^{q})^{2\alpha-1}\\ \times
(\xi_{-}/2^{1/q}(1+l)/2)^{q-p-1}dl.
\end{eqnarray*}
The same approach is applied to all integrals in the previous 
section. Note that in contrast to multi-domain approaches for 
differential equations, no gluing of the different domains is needed 
here since the antiderivative is a continuous map. Thus the integrals 
can simply be added after being computed with the Clenshaw-Curtis 
algorithm.

\subsection{Example}
It is instructive to test the code for an explicitly known example. 
We consider
\begin{equation}
	D^{1/2}\frac{1}{1+x^{2}}= \frac{\sqrt{\pi}}{4}\left(
	\frac{1}{(1-ix)^{3/2}}+\frac{1}{(1+ix)^{3/2}}\right),
	\label{exact}
\end{equation}
i.e., a function that has a stronger fall-off at infinity than 
required by our approach 
($\lim_{x\to\infty}\frac{|x|^{3/2}}{1+x^{2}}=0$). The fractional 
derivative can be computed either directly from the integral 
(\ref{fracint}) or the 
Fourier transform. 

\begin{figure}[htb!]
 \includegraphics[width=0.49\textwidth]{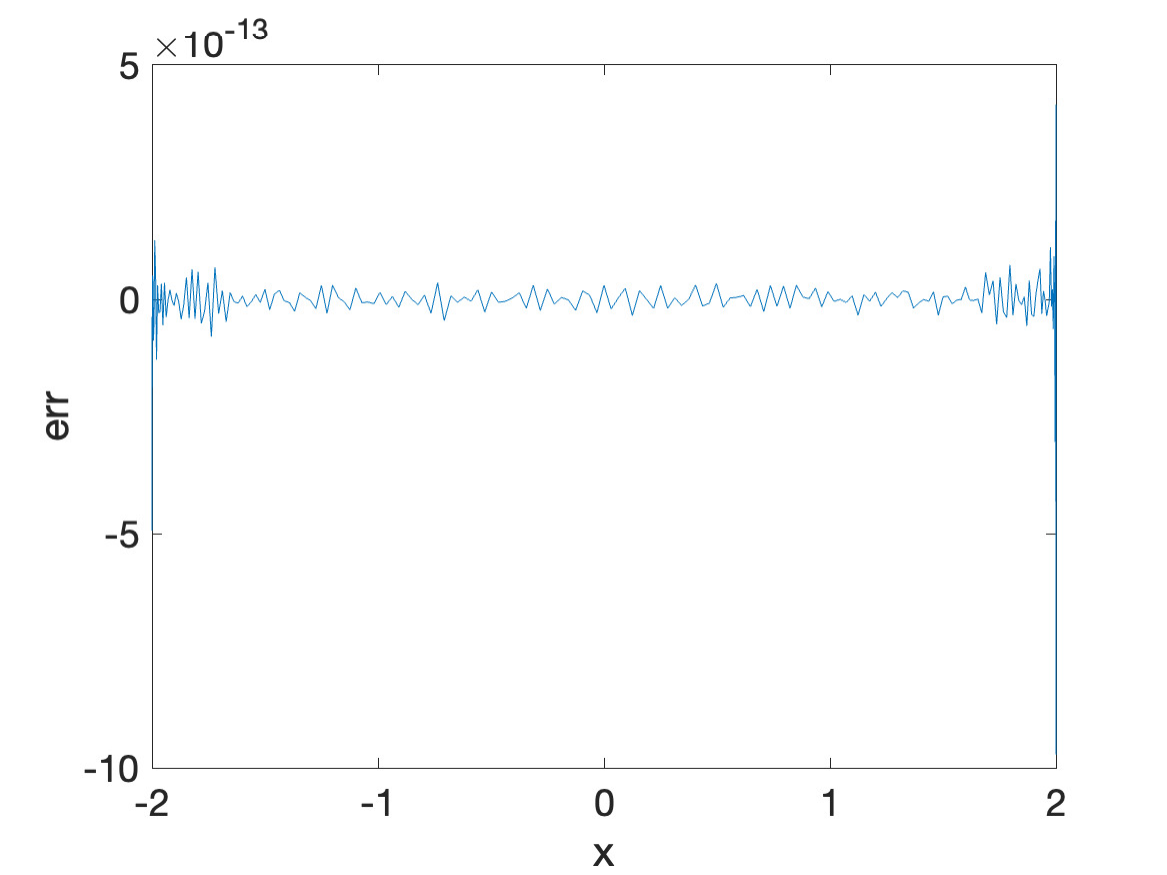}
 \includegraphics[width=0.49\textwidth]{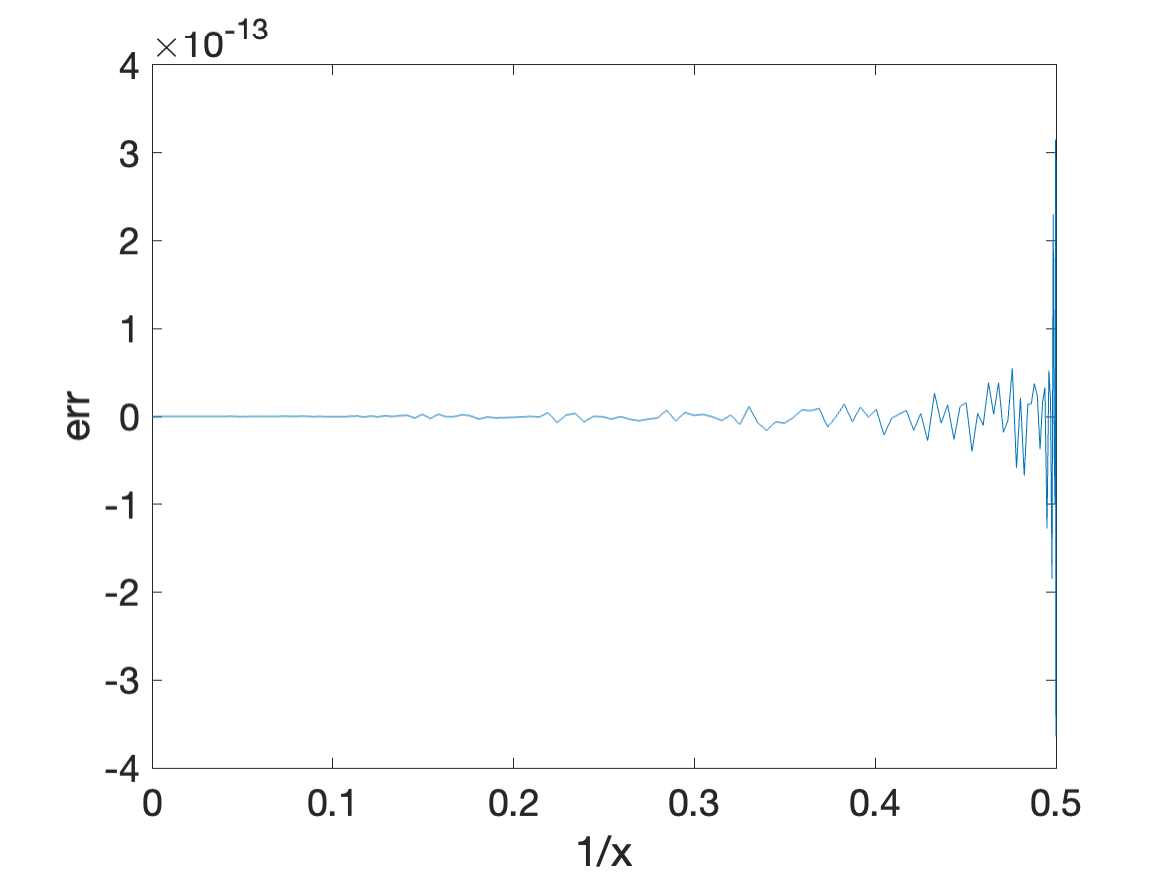}
 \caption{Difference of the numerical solution for the fractional 
 derivative (\ref{exact}) and its exact expression, on the left for 
 domain II, on the right for domain III.}
 \label{figexample}
\end{figure}

To treat this example numerically, we choose $b=-a=2$ and $N=200$ 
collocation points for each integral. The difference between the 
numerical solution and the exact solution can be seen in 
Fig.~\ref{figexample} for domain II on the left and domain III on the 
right. It is of the order of $10^{-13}$ as expected. Since the 
function in (\ref{exact}) is symmetric, this is also the case for the 
fractional derivative. In fact $I^{+}$ in domain III coincides with 
$-I^{-}$ in domain I and vice versa up the order of $10^{-16}$. The 
numerical differentiation introduces an error of the order of 
$10^{-12}$ near infinity. Thus the error in domain I is not shown 
since it is similar to the one in domain III.  

\subsection{Comparison with DFT}
Probably the most used spectral approach to compute fractional derivatives 
(\ref{fourier}) consists in approximating the expression in Fourier 
space by a DFT. We will now study to what extent this introduces 
numerical errors. 

\begin{figure}[htb!]
 \includegraphics[width=0.49\textwidth]{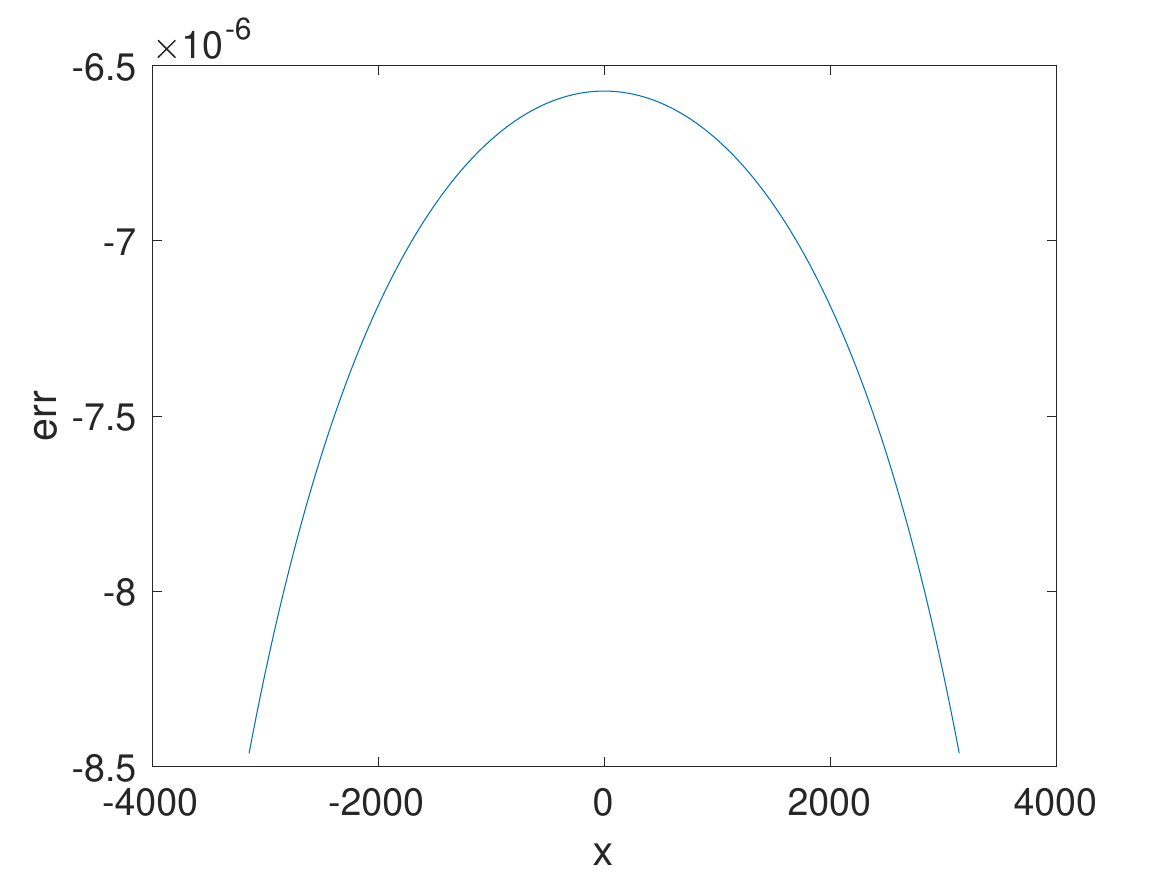}
 \includegraphics[width=0.49\textwidth]{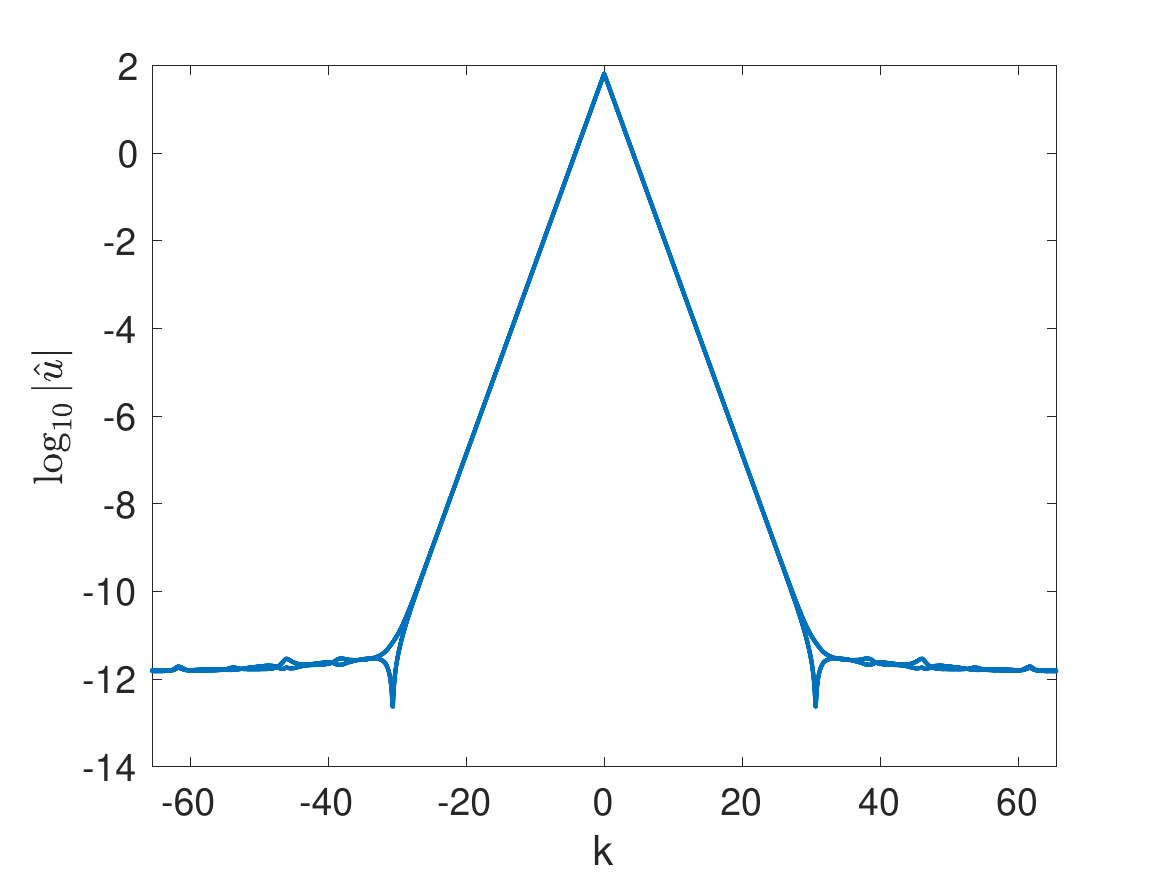}
 \caption{Difference of the numerical solution for the fractional 
 derivative (\ref{exact}) with an FFT approach 
 and its exact expression on the left, and the DFT coefficients for 
 the fractional derivative on the right.}
 \label{errFFT}
\end{figure}

To treat the exactly known example (\ref{exact}), we work for 
$x\in L[-\pi,\pi]$ with $L=1000$ with $N_{FFT}=2^{17}\approx 1.3*10^{5}$ DFT modes. Here we 
apply the standard DFT discretisation,
\begin{equation}
	x_{n} = -\pi L +nh,\quad n = 1,\ldots,N_{FFT},\quad h = 2\pi 
	L/N_{FFT},
	\label{FFT}
\end{equation}
and $k = (-N_{FFT}/2+1,\ldots,N_{FFT}/2)/L$. As can be seen in 
Fig.~\ref{errFFT} on the left, the 
difference between numerical and exact solution denoted by err, 
is of the order $10^{-6}$. 
Note that this error is not related to the decrease of the DFT 
coefficients $\hat{u}$ shown on the right of the same figure 
in a semi-logarithmic plot. They 
decrease to the order of $10^{-10}$. In \cite{KPS}, the difference 
between an FFT based approach to the Lorentz profile (\ref{exact}) and 
a Fourier approach on the whole real line was shown to be of the 
order of the highest DFT coefficients, here $10^{-10}$. The reason 
why the error is higher for fractional derivatives is that there are 
Abelian integrals from $-\infty$ to $-L\pi$ and from $L\pi$ to 
$\infty$ systematically neglected by an approximation on the torus. 
This is reflected by the non-analytic term $|k|^{\alpha}$ in the 
definition of the fractional derivatives via Fourier transforms. 
As indicated by the DFT coefficients in Fig.~\ref{errFFT}, the 
numerical error in approximating the fractional derivative is not due 
to a lack of resolution in the DFT computation. A dependence of the 
numerical error can be seen in table \ref{table1} on the left. For values of 
$N_{FFT}=2^{16}, 2^{17}, 2^{18}, 2^{19}$, the maximal 
difference (err) between DFT computation and the exact expression is 
as for $N_{FFT}=2^{15}$. The dependence on the size of the 
torus is shown in table \ref{table1}  on the right where we have used $N_{FFT}=2^{19}$. 
\begin{table}[tbp]
	\centering
		\begin{tabular}{|c|c|}
		\hline
		$N_{FFT}$ & err   \\
		\hline
		$2^{12}$ & $8.25*10^{-2}$   \\
		\hline
		$2^{13}$ & $8.43*10^{-3}$   \\
		\hline
		$2^{14}$ & $1.04*10^{-4}$   \\
		\hline
		$2^{15}$ & $8.46*10^{-6}$   \\
		\hline
	\end{tabular}
	\qquad
	\begin{tabular}{|c|c|}
		\hline
		$L$ & err   \\
		\hline
		$10^{2}$ & $2.68*10^{-4}$   \\
		\hline
		$10^{3}$ & $8.64*10^{-6}$   \\
		\hline
		$10^{4}$ & $2.68*10^{-7}$   \\
		\hline
		$10^{5}$ & $4.38*10^{-2}$   \\
		\hline
	\end{tabular}
	\caption{$L^{\infty}$ norm of the difference of the numerical solution for the fractional 
 derivative (\ref{exact}) with an FFT approach 
 and its exact expression for $N_{FFT}=2^{19}\approx 5.2*10^{5}$ and different sizes of 
 the torus.}
		\label{table1}
\end{table}

\begin{figure}[htb!]
 \includegraphics[width=0.49\textwidth]{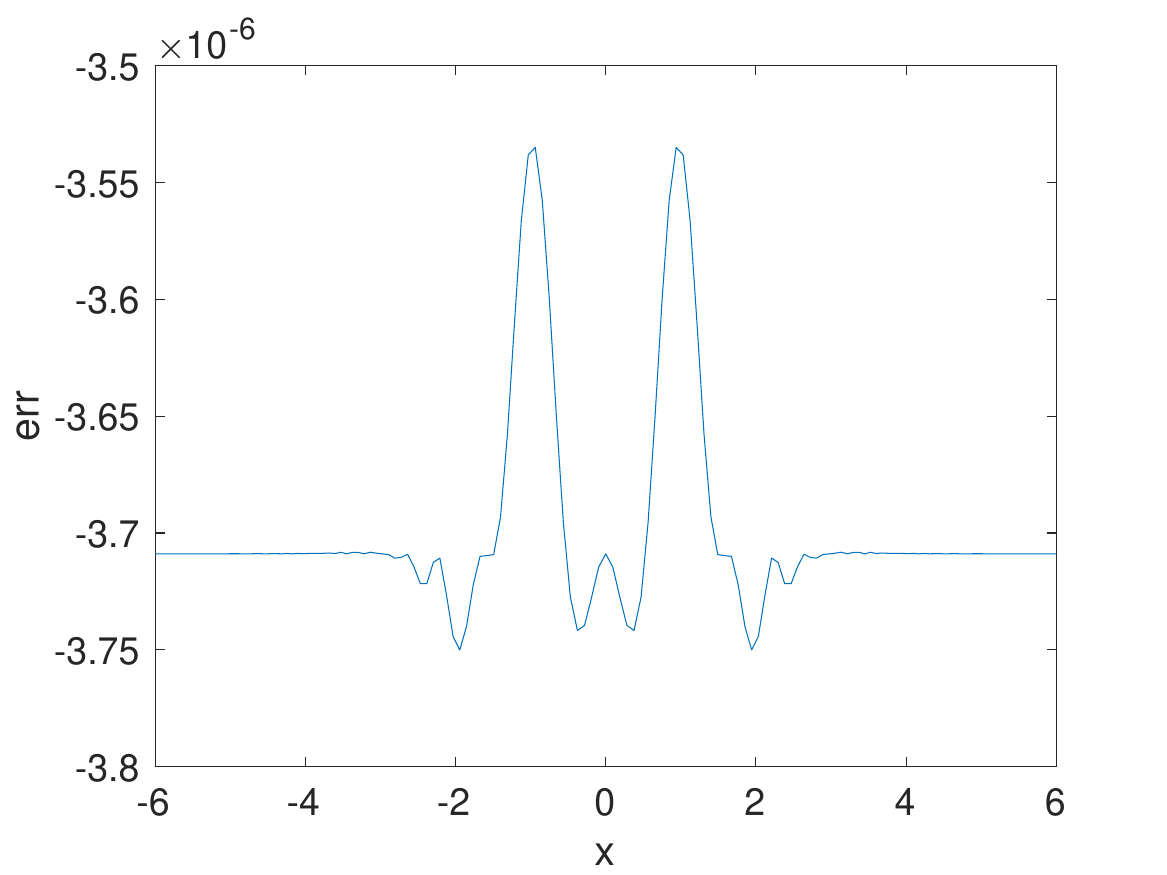}
 \includegraphics[width=0.49\textwidth]{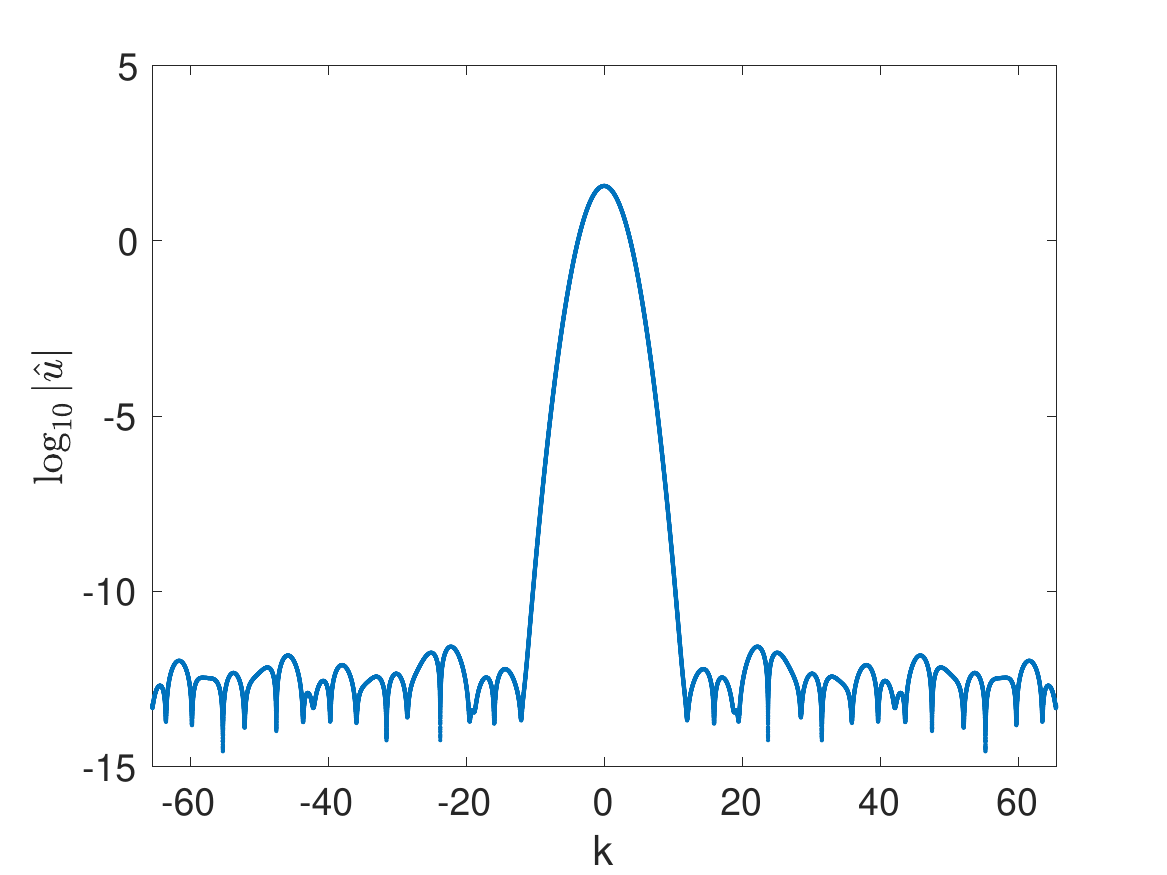}
 \caption{Difference of the numerical solution for the fractional 
 derivative of order $1/2$ of a Gaussian with an FFT approach 
 and the multi-domain approach on the left, and the DFT coefficients for 
 the Gaussian on the right.}
 \label{errFFTgauss}
\end{figure}

The fact that the applicability of the DFT in this context is not 
limited by the decrease of the function $u$ for $x\to\infty$ can be 
further illustrated by considering fractional derivatives of Schwartz 
class functions, for instance $D^{1/2}\exp(-x^{2})$. For 
$x=0$, this derivative is known to be given by 
$4^{3/4}\Gamma(3/4)/(2\sqrt{\pi)})$, where $\Gamma$ is the Gamma 
function. If we use the multi-domain approach with $b=-a=6$ and 
$N=200$ (in domains I and III the Gaussian is smaller than $10^{-16}$ 
and thus vanishes to numerical precision; this implies that the 
related integrals to $u^{I}$, $u^{III}$ can be simply ignored),
the difference of exact (computed with the Matlab implementation of 
the $\Gamma$ function) and numerical solution for $x=0$ is 
$8.55*10^{-15}$, as expected of the order of machine precision. In 
contrast the FFT approach with $N_{FFT}=2^{17}$ DFT modes for $x\in 
1000[-\pi,\pi]$, leads to a difference of numerical and exact 
solution of $3.71*10^{-6}$ as before. This illustrates that the 
numerical error is not due to the approximation of the function $u$ 
via DFT
(which is here of the order of machine precision as can be seen from 
the DFT coefficients on the right of Fig.~\ref{errFFTgauss}) but to the singular 
symbol $|k|^{\alpha}$ in the inverse Fourier transform. The 
difference between the fractional derivative computed with the 
multi-domain approach and the DFT approach can be seen in 
Fig.~\ref{errFFTgauss} on the left to be of the order of $10^{-6}$ as 
in Fig.~\ref{errFFT}. 

\begin{figure}[htb!]
 \includegraphics[width=0.49\textwidth]{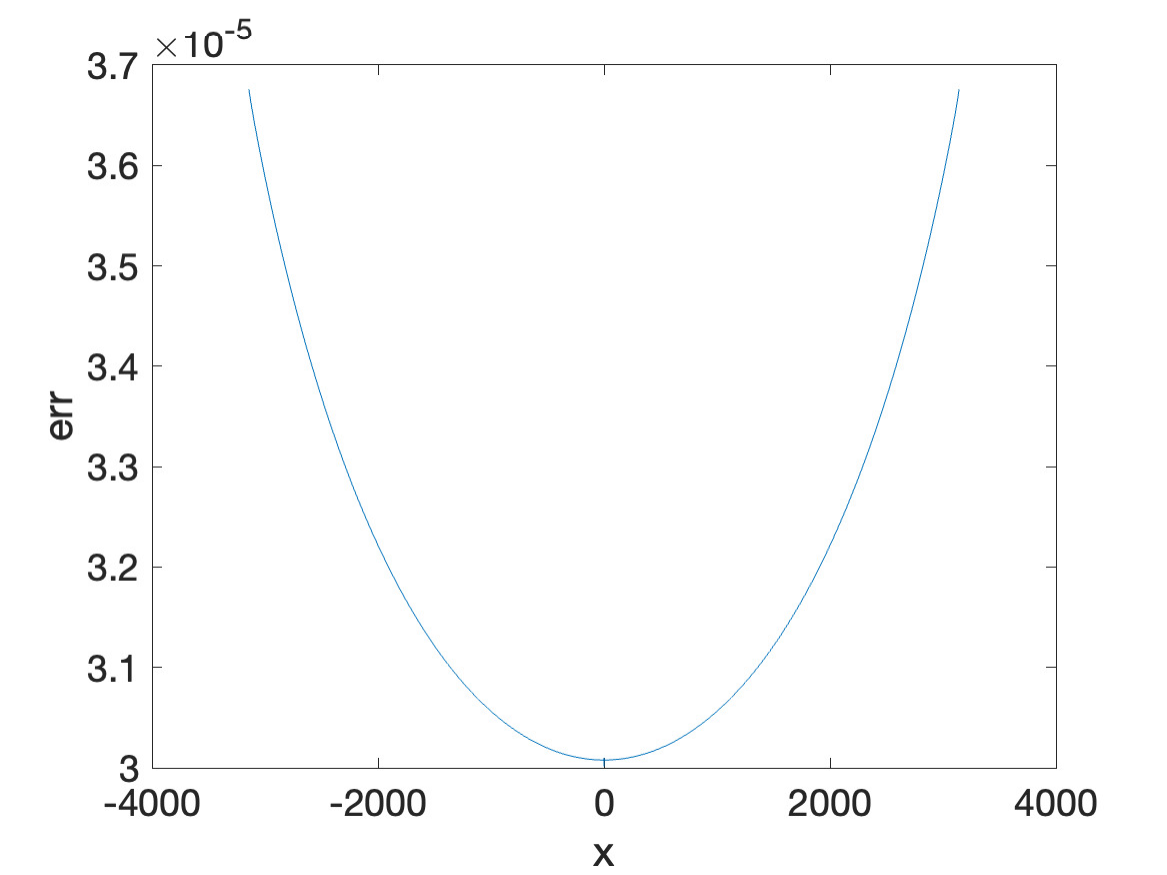}
 \includegraphics[width=0.49\textwidth]{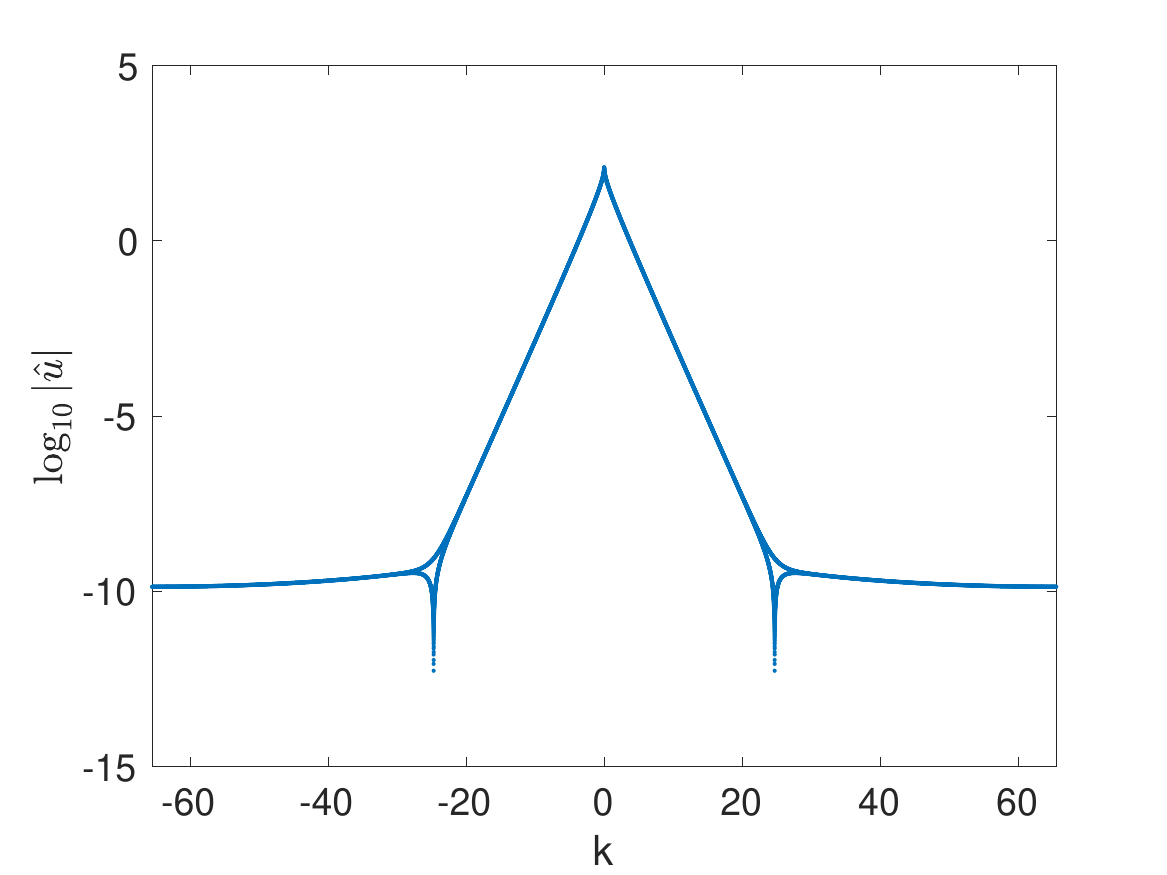}
 \caption{Difference of the numerical solution for the fractional 
 derivative (\ref{example}) for $\alpha=2/5$ with an FFT approach 
 and the multi-domain approach on the left, and the DFT coefficients for 
 the fractional derivative on the right.}
 \label{errexample}
\end{figure}

It was shown in \cite{KPS} that the error increases if the fall-off 
of the considered function is slower, for instance as in the case of 
the fractional KdV solitary waves. If we consider an example of this 
kind, for instance
\begin{equation}
	u(x)= \frac{1}{(1+x^{2})^{(1+\alpha)/2}}
	\label{example},
\end{equation}
the discrepancy between an FFT based approach and the multi-domain 
approach in this paper can be even larger. We are not aware of an 
exact expression in such a case and consider therefore the 
difference between the two approaches. With the same numerical 
parameters as before, we get for $\alpha=2/5$ the difference between multi-domain and 
FFT approach on the left of Fig.~\ref{errexample}. It is of the order 
of $10^{-5}$. As the DFT coefficients on the right of the same 
figures show, the error is considerably larger than the magnitude of 
the DFT coefficients with the largest index. Thus the error does not 
change much if a higher resolution or a larger torus is used for the 
DFT. 

The optimal implementation of the DFT in particular 
in Matlab means that the computation of a fractional derivative with 
FFT techniques is much more rapid. But as the examples above show, 
the achievable accuracy is limited to roughly $10^{-6}$. The 
multi-domain approach on the other hand is more complicated and takes 
more computation time, but can reach essentially machine precision.

\section{Solitary waves of fractional KdV equations}
In this section, we construct solitary wave solutions to the 
fractional KdV equation with the methods of the previous section. 
Equation (\ref{fKdVsol}) is solved with the multi-domain approach and 
a Newton-Krylov iteration. The results are 
compared to what was obtained with an FFT approach in \cite{KS15}. 
For an application of the 
Petviashvili method with FFT in this context, the reader is referred 
to \cite{Dur}. Note that is was shown in \cite{LPS} that such 
solitary waves only exist for $\alpha>1/3$ for fKdV.

\subsection{Fractional Korteweg-de Vries equation}
We refer the reader to \cite{KSbook} for a review on the fractional 
KdV equation (\ref{fKdV}).
Here we will just summarize important features needed in the 
following. For $\alpha=2$, this is the standard KdV equation, for 
$\alpha=1$ the Benjamin-Ono equation (for the case $\alpha=1$ we do not 
consider here, the fractional 
derivative $D^{\alpha}$ becomes the Hilbert transform).

The following quantities are formally conserved by the fractional KdV 
flow: the mass
\begin{equation} \label{M}
M(u)=\int_{\mathbb R}u^2(x,t)dx,
\end{equation}
and the Hamiltonian
\begin{equation} \label{H}
H(u)=\int_{\mathbb R}\big( \frac{1}{2} |D^{\frac{\alpha}2}u(x,t)|^2-\frac{1}{6}u^3(x,t)\big) dx.
\end{equation}

The fractional KdV  equation (\ref{fKdV}) is invariant under the scaling transformation 
\begin{equation} \label{scaling}
u_{\lambda} (x,t)=\lambda^{\alpha}u(\lambda x,\lambda^{\alpha+1}t),
\end{equation}
for any positive $\lambda\in\mathbb{R}$. The $L^{2}$ norm is 
invariant under this scaling for $\alpha=\frac{1}{2}$, this is  called the 
$L^2$-critical case. The energy is invariant for $\alpha=1/3$, the 
energy-critical case.

Solitary wave solutions of (\ref{fKdV}) are localized traveling wave 
solutions, i.e., solutions of the form $u(x,t)=Q_{c}(x-ct)$ where $c>0$ 
is a constant velocity. With the fall-off condition at infinity, 
the solitary waves of fractional KdV are given after  integration 
by (\ref{fKdVsol}). There are no solitary waves of finite mass and energy 
for $\alpha\leq 1/3$. It is sufficient to concentrate on solitary 
waves for $c=1$ since $Q_{c}$ is related to $Q:=Q_{1}$ via the 
scaling 
\begin{equation}
	Q_{c}(z) = cQ(zc^{1/\alpha}).
	\label{Qc}
\end{equation}
For $\alpha=1$, the solitary wave is known explicitly, 
$Q(z)=4/(1+z^{2})$, the soliton of the Benjamin-Ono equation. For general $\alpha\in ]0,1[$, it was shown in 
\cite{FL} that the fall-off behavior of $Q(z)$ for $|z|\to\infty$ is 
proportional to $1/|z|^{1+\alpha}$. 

\subsection{Numerical construction of solitary waves}
Equation (\ref{fKdVsol}) was solved numerically in \cite{KS15} by 
considering it on a torus, $x\in L[-\pi,\pi]$ with $L=100$ and 
discretising it with the standard FFT discretisation (\ref{FFT}). The 
resulting set of nonlinear equations for $Q$ of the form $F(Q)=0$ was then solved with a 
standard Newton iteration, 
$$Q^{(n+1)}=Q^{(n)}-\left((\mathrm{Jac}(F)|_{Q^{(n)}}\right)^{-1}F(Q^{(n)}),$$
where $\mathrm{Jac}(F)$ is the Jacobian of $F$. 
In an abuse of notation, we have denoted the (discrete) iterates $Q^{(n)}$, 
$n=0,1,\ldots$ with the same symbol as the continuous variables. The 
action of the inverse of the Jacobian on the vector $F$ is computed iteratively with 
the Krylov subspace technique GMRES\footnote{The basic idea of 
GMRES is that the solution of the equation $Ax=b$, $x,b\in 
\mathbb{R}^{N}$, $A$ an invertible $N\times N$ matrix, is expressed in terms of 
a linear combination of
the vectors $A^{n}b$, $n=0,1,\ldots$, see \cite{GMRES} for details. 
The advantage of this approach is that only the action of $A$ on $b$ 
needs to be known, not $A$ itself.}. Note that GMRES can fail to 
converge to the prescribed tolerance (here $10^{-13}$), but this will not necessarily imply that the Newton 
iteration will not. It is well known that a Newton iteration with an 
approximate Jacobian may still converge, but in general only linearly 
and not quadratically as with the correct Jacobian. 

Since the convergence of the Newton iteration is local, the choice of 
the initial iterate $Q^{(0)}$ is important. In \cite{KS15}, we 
applied a tracing technique: we started with the explicitly known 
Benjamin-Ono soliton as the initial iterate for $\alpha=0.9$, then 
used the resulting solution for $\alpha=0.8$ and so on. The closer 
one gets to the energy-critical case, the smaller these steps have to 
be chosen. In \cite{KS15} we used $N_{FFT}=2^{16}$ Fourier modes and 
could reach $\alpha=0.45$ as the smallest value. 

In the present paper we use the multi-domain approach detailed in the 
previous sections to discretise equation (\ref{fKdVsol}). This leads 
once more to a system of nonlinear equations to be solved with a 
Newton-Krylov method. Note that we work in domains I and III with the 
functions $u^{I}$ and $u^{III}$ defined in (\ref{cases}). We use the solution obtained with an FFT 
approach as the initial iterate 
(for $|x|>L\pi$, the functions $u^{I}$, $u^{III}$ are replaced with the constant value at $\pm L\pi$ respectively, thus insuring the appropriate fall of at $\pm\infty$).
 Note that the 
solitary wave can be chosen to be symmetric with respect to $x\mapsto 
-x$. We do not use this symmetry in the numerical approach. 

\subsection{fKdV solitary waves for moderate values of $\alpha$}
\begin{figure}[htb!]
 \includegraphics[width=0.7\textwidth]{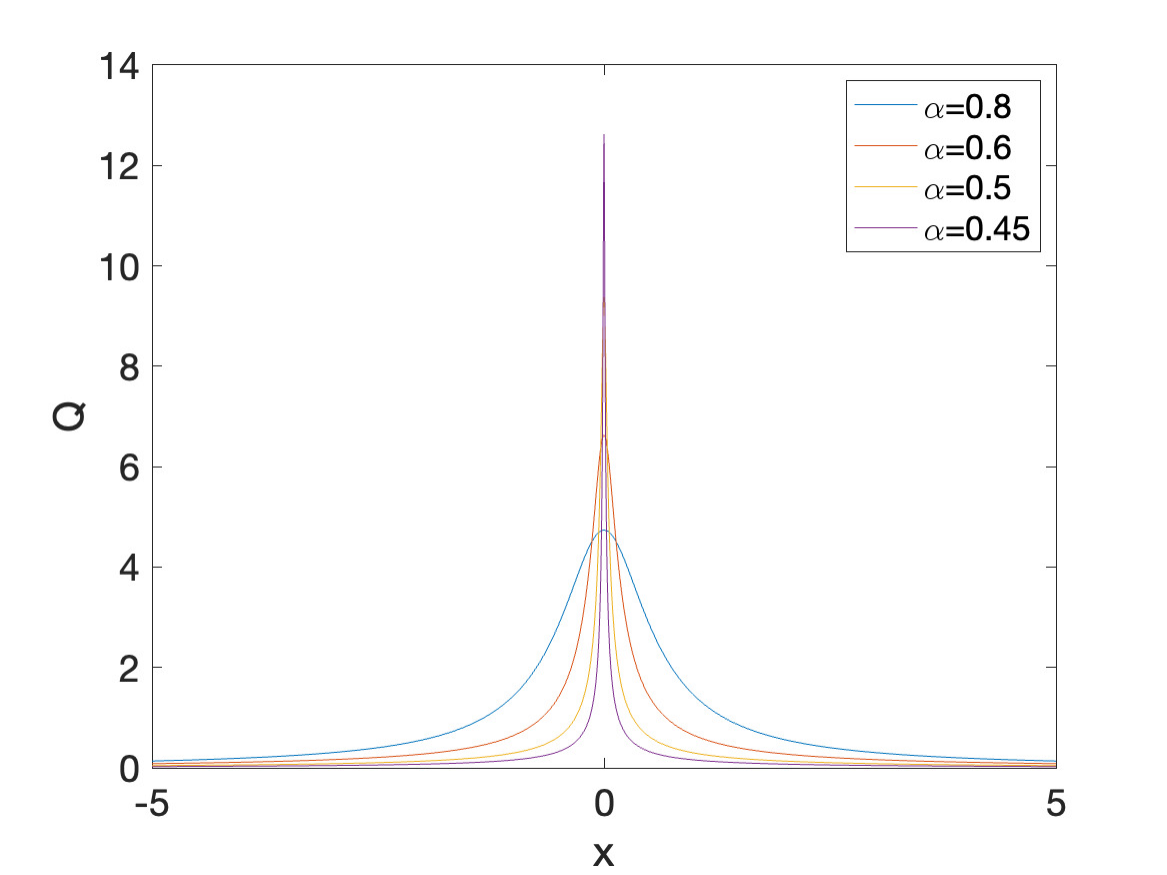}
 \caption{Solitary wave of the fKdV equation for $c=1$ and several 
 values of 
 $\alpha$.}
 \label{fKdVsolalpha}
\end{figure}

We discuss first the case $\alpha=0.8$ with $N=200$, $b=-a=1$, and 
$\delta=10^{-2}$. Starting with the result of 
the FFT iteration (interpolated to the Chebyshev grid), the Newton 
iteration in the multi-domain case converges 
in 4 iterations to a residual of the order of $10^{-12}$. The 
solution is shown for several values of $\alpha$ in 
Fig.~\ref{fKdVsolalpha}. For $\alpha=0.8$ (blue), it is still very 
close to the Lorentz profile for the Benjamin-Ono soliton 
($\alpha=1$). 

\begin{figure}[htb!]
 \includegraphics[width=0.49\textwidth]{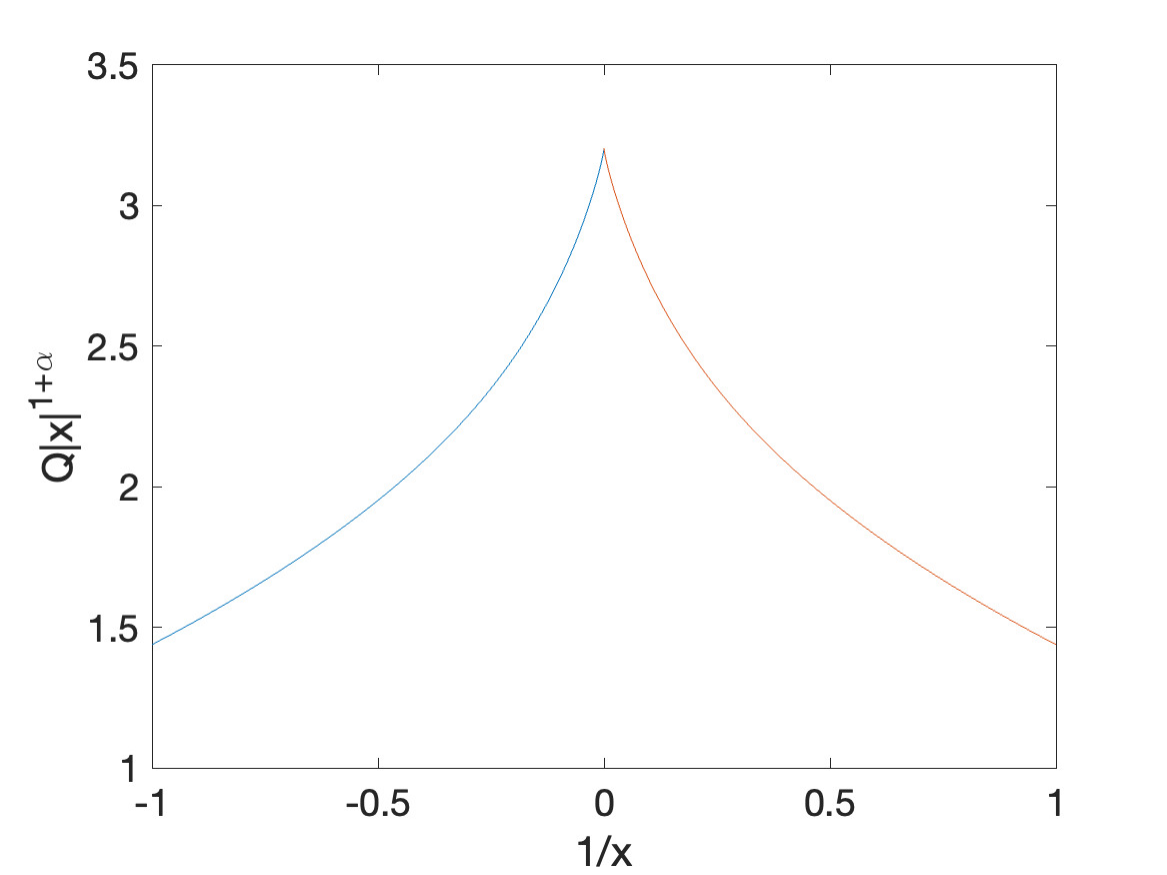}
 \includegraphics[width=0.49\textwidth]{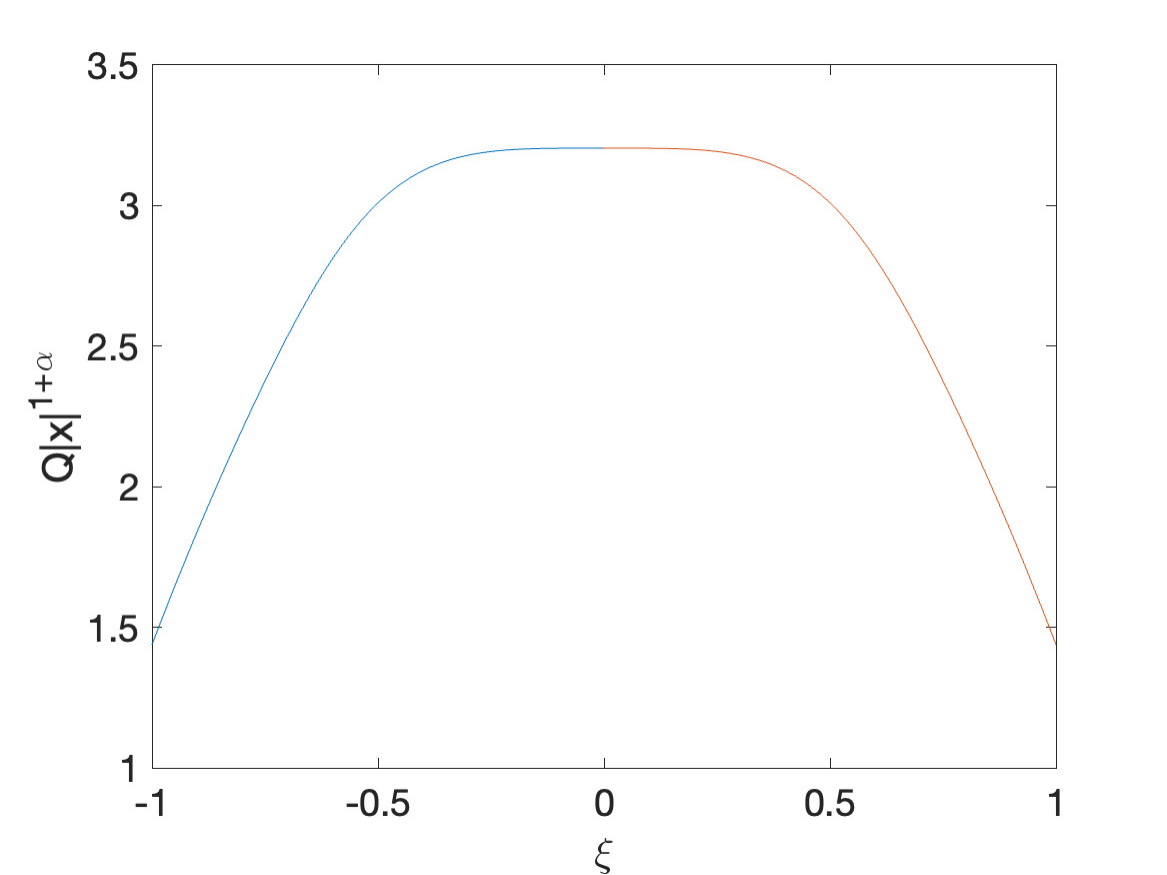}
 \caption{Solitary wave of the fKdV equation for $c=1$ and 
 $\alpha=0.8$ near infinity, on the left in dependence of $1/x$, on 
 the right in dependence of $\xi_{+}$ and $-\xi_{-}$.}
 \label{fKdVsol08inf}
\end{figure}

The solution near infinity can be seen in Fig.~\ref{fKdVsol08inf}. We 
show the functions $Q^{I}$ and $Q^{III}$ in domains I and III 
respectively, i.e., $Q|x|^{1+\alpha}$. On the left of the figure, we 
show these functions in dependence of $1/x$. They 
tend to a constant nonzero value at infinity which is in our approach 
just a point on the grid. It can also be seen that the functions have 
a cusp at infinity as expected from the fact that the fractional 
derivatives are defined on a Riemann surface. This means that the 
slopes diverge at infinity which makes the application of a spectral 
method in $1/x$ inefficient (thus preventing us, for example, to use 
a single infinite domain in $1/x$). However, the solutions are as expected 
smooth in an adapted local parameter $\xi_{\pm}=1/(\pm x)^{1/q}$ at 
$\xi_{\pm}=0$, i.e., near 
infinity on the Riemann surface. This can be seen on the right of the 
same figure (for this representation, the solution is shown in 
dependence of $-\xi_{-}$ and $\xi_{+}$). 

\begin{figure}[htb!]
 \includegraphics[width=0.49\textwidth]{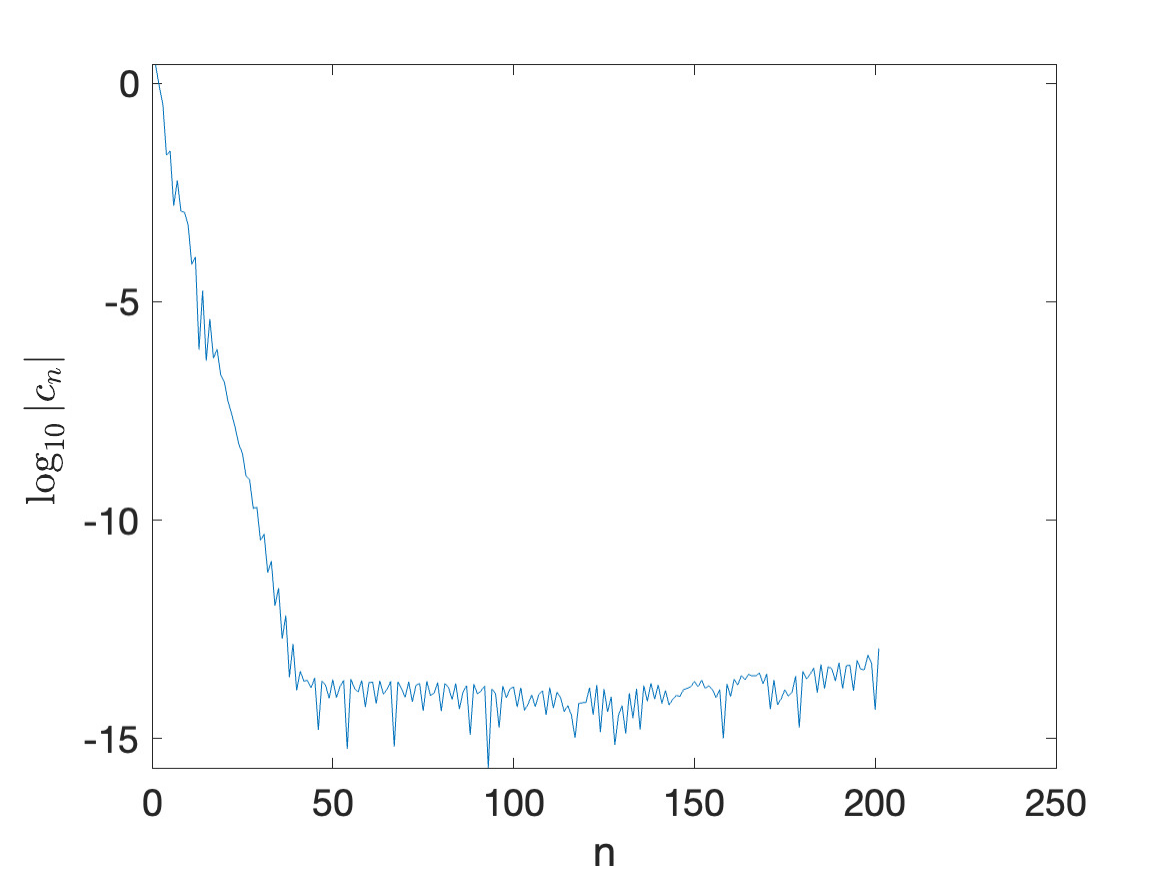}
 \includegraphics[width=0.49\textwidth]{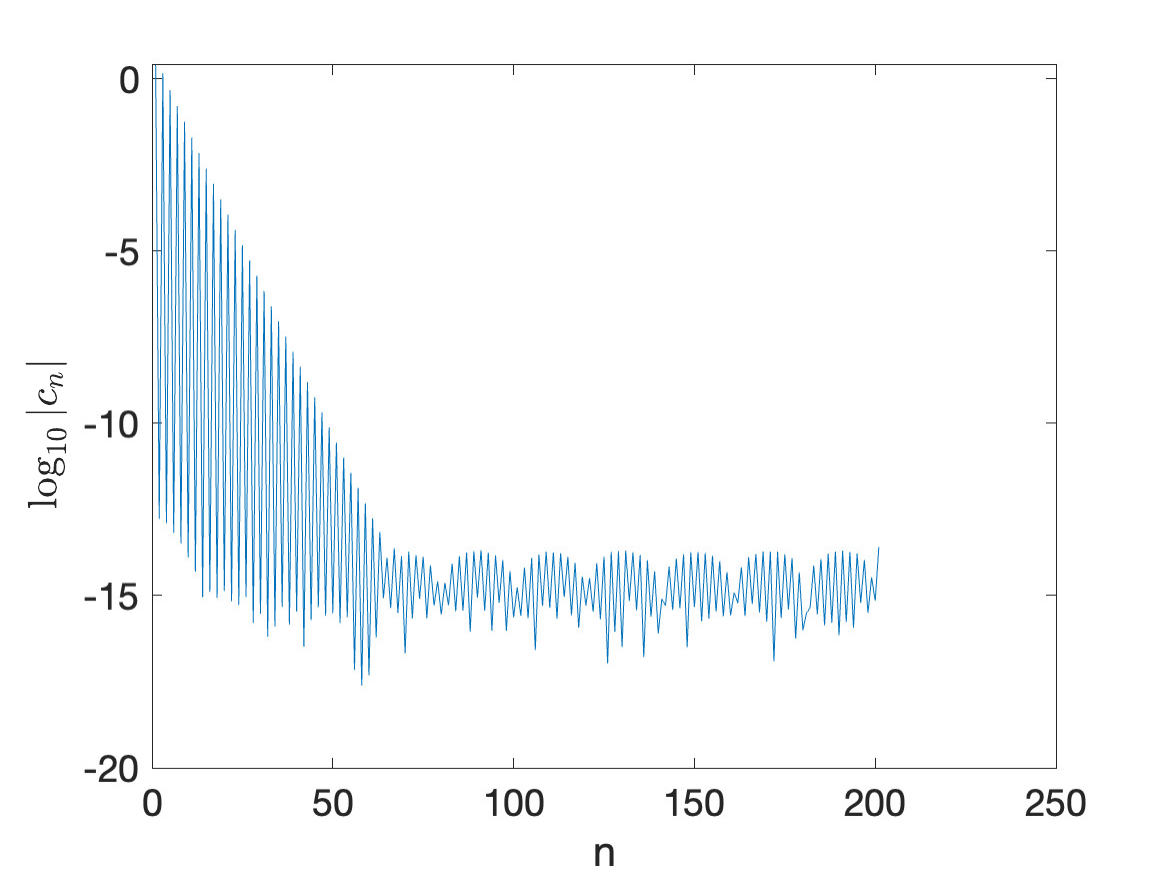}
 \caption{Chebyshev coefficients for the solitary wave of the fKdV equation for $c=1$ and 
 $\alpha=0.8$, on the left in domain I, on 
 the right in domain II.}
 \label{fKdVsol08cheb}
\end{figure}

\begin{figure}[htb!]
 \includegraphics[width=0.49\textwidth]{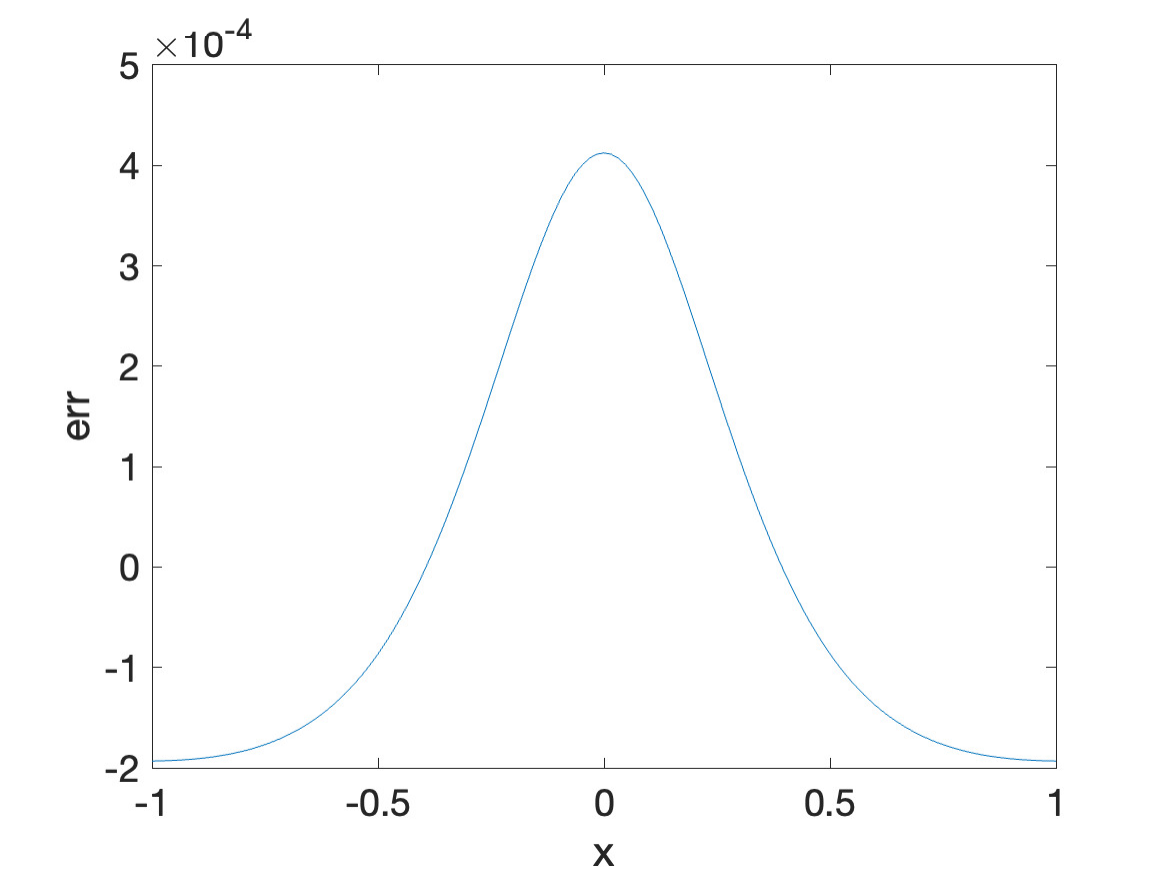}
 \includegraphics[width=0.49\textwidth]{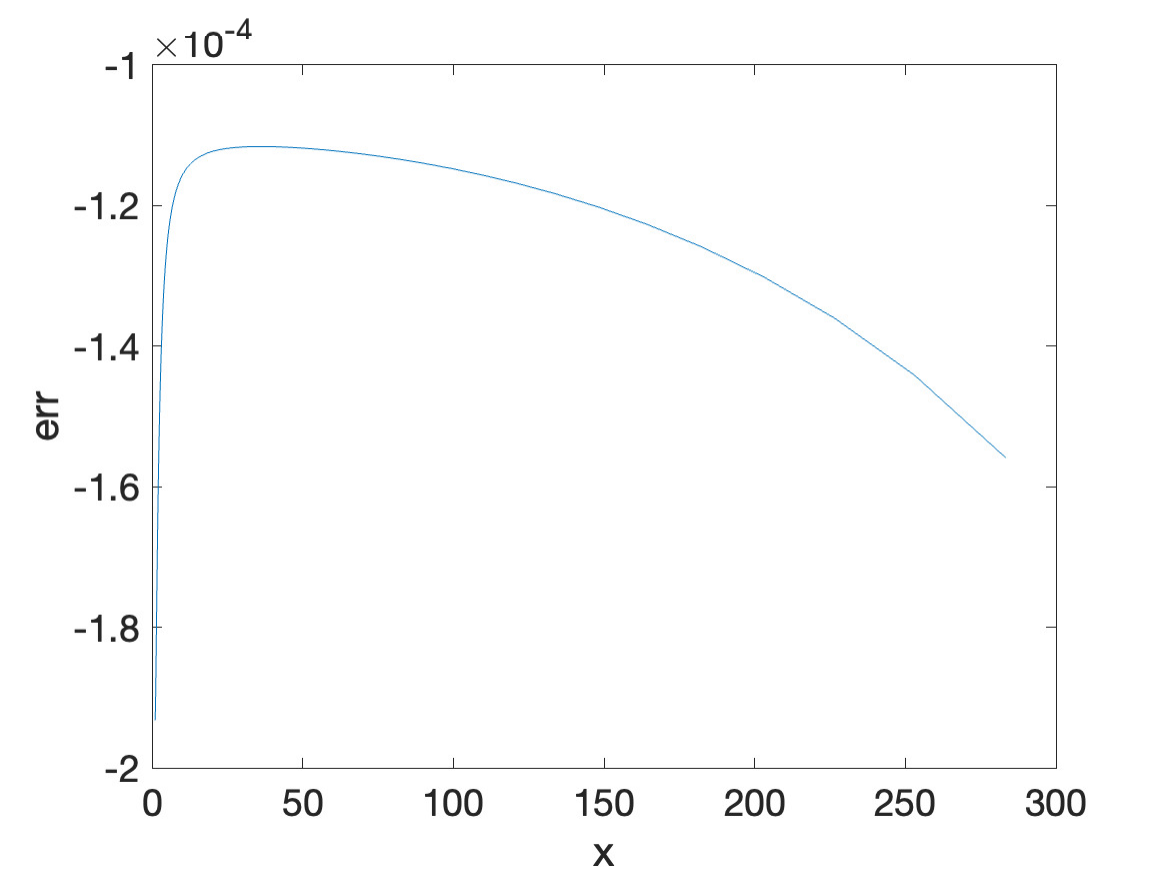}
 \caption{Difference between the solitary wave of the fKdV equation for $c=1$ and 
 $\alpha=0.8$ constructed with an FFT and a multi-domain approach for 
 various values of $x$.}
 \label{fKdVsol08FFT}
\end{figure}

A spectral approximation of $Q$ in dependence of $\xi_{\pm}$ is 
therefore very efficient as shown by the Chebyshev coefficients in 
Fig.~\ref{fKdVsol08cheb}. It can be seen that $N\sim 70$ would be 
sufficient to fully resolve the solitary wave in coefficient space 
(the highest coefficients are of the order of machine precision). 
However, we need values of $N>160$ in the computation of the 
fractional integrals to reach machine precision even if this is not 
needed to resolve the final result in the space of Chebyshev 
coefficients. Note that we only show the coefficients in domain I and 
II since they are the same in domain III and I up to numerical errors 
for symmetry reasons. This also explains why every other coefficient 
in domain II vanishes with numerical precision. A more efficient 
approach would be to use two finite domains ($[a,0]$ and $[0,b]$) or 
to explicitly implement the symmetry of $Q$, $Q(-x)=Q(x)$. 
It is interesting to note the difference between the solution via the multi-domain
approach and the solitary wave constructed with an FFT approach. This difference is shown in 
Fig.~\ref{fKdVsol08FFT}. It is of the order of $10^{-4}$. Note that 
the comparison of fractional derivatives computed with FFT and with 
the multi-domain approach was of the order of $10^{-6}$ in the 
previous section, mainly because the restriction to a finite torus 
neglects the fractional integrals from $\pm \pi L$ to infinity. Here 
the effect is amplified by it appearing within an iteration. 

\subsection{fKdV solitary waves for values of $\alpha$ closer to the 
critical value $\alpha_{c}=1/3$}

\begin{figure}[htb!]
 \includegraphics[width=0.49\textwidth]{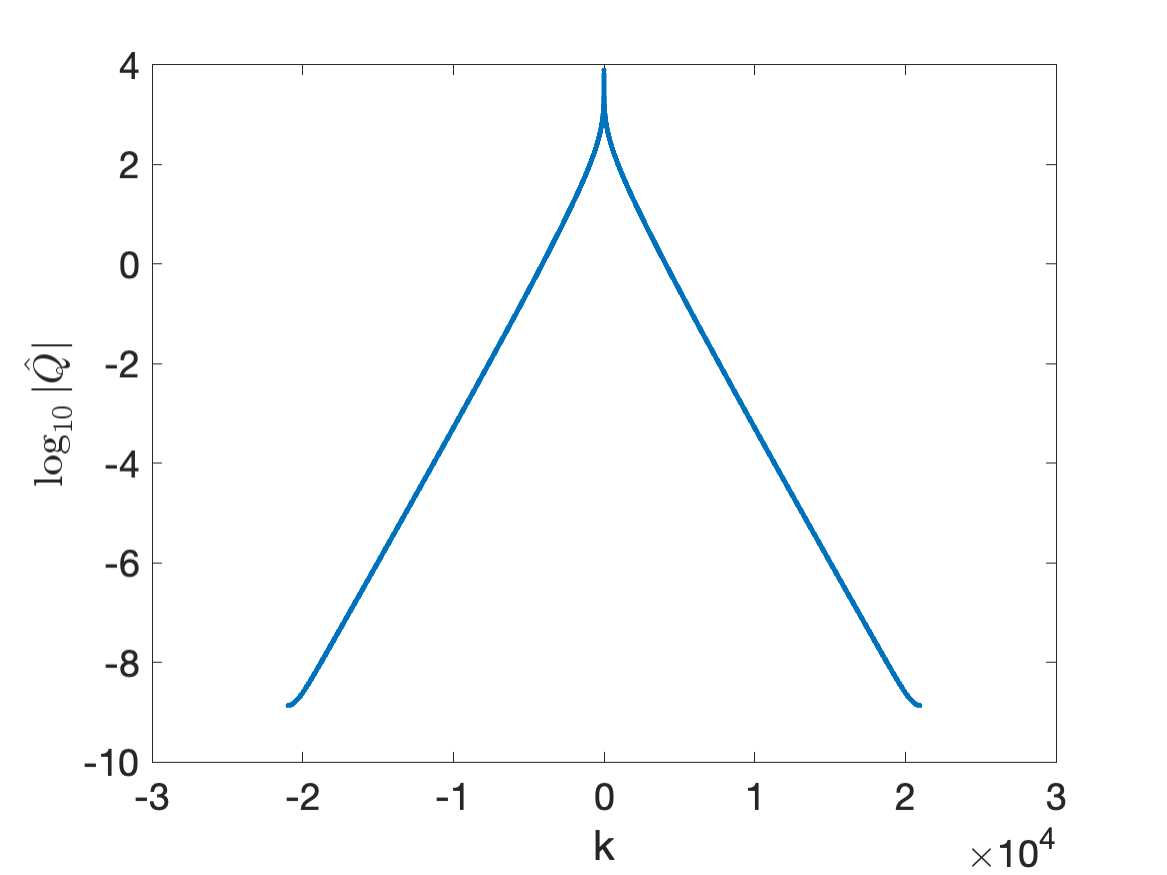}
 \includegraphics[width=0.49\textwidth]{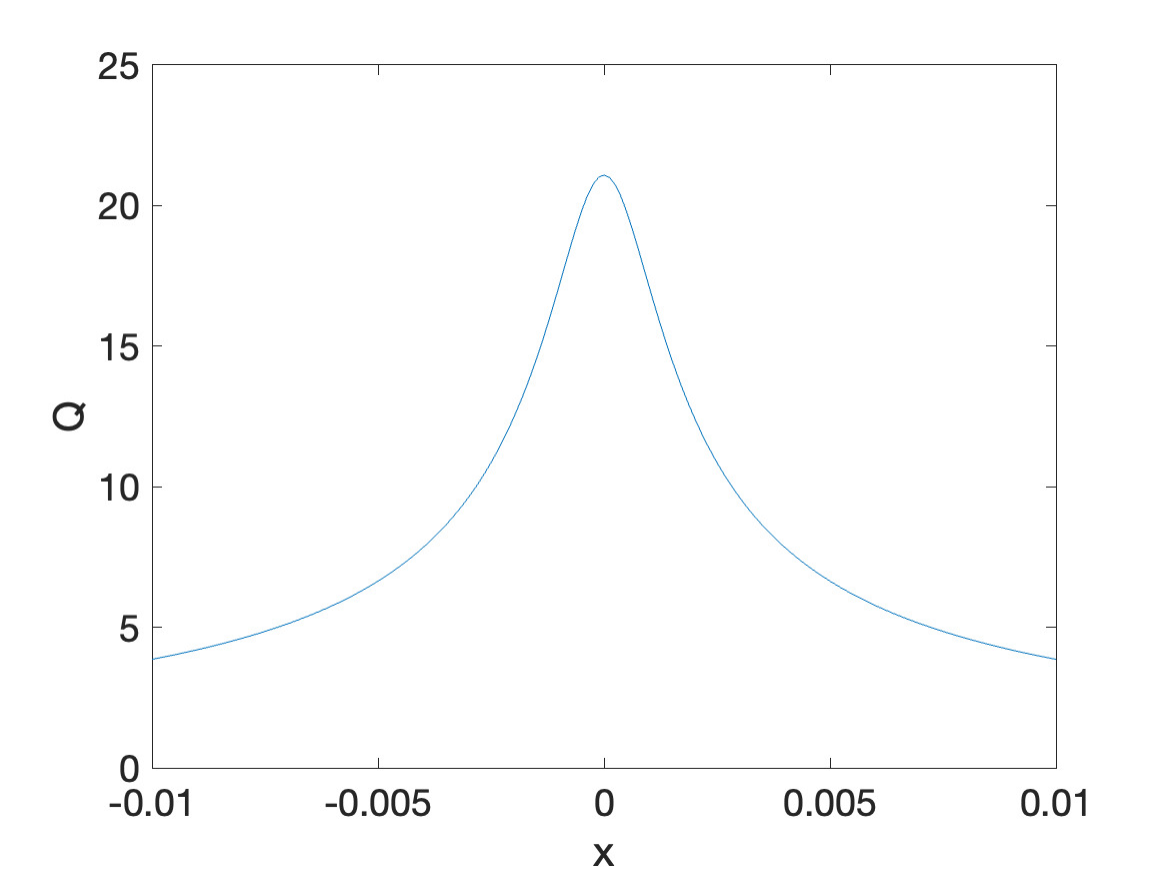}
 \caption{Solitary wave of the fKdV equation for $c=1$ and 
 $\alpha=0.4$: on the left the DFT coefficients for the solution 
 constructed with an FFT approach, on the right the solution to the 
 multi-domain approach.}
 \label{fKdVsol04FFT}
\end{figure}

In \cite{KS15}, solitary waves were constructed for values of 
$\alpha\geq0.45$. Examples can be seen in Fig.~\ref{fKdVsolalpha}. 
The smaller $\alpha$, the more the solution is peaked at $x=0$, but 
with shrinking support and slower fall-off to infinity. 
With $N_{FFT}=2^{22}\approx 4.2*10^{6}$, one can reach values closer to the 
critical value $\alpha_{c}=1/3$ below which no smooth solitary waves 
exist. The DFT coefficients of the solitary wave for $\alpha=0.4$ are 
shown in Fig.~\ref{fKdVsol04FFT} on the left. They decrease to the 
order of $10^{-8}$. We could not reach even smaller values of 
$\alpha$ with this approach. 

\begin{figure}[htb!]
 \includegraphics[width=0.49\textwidth]{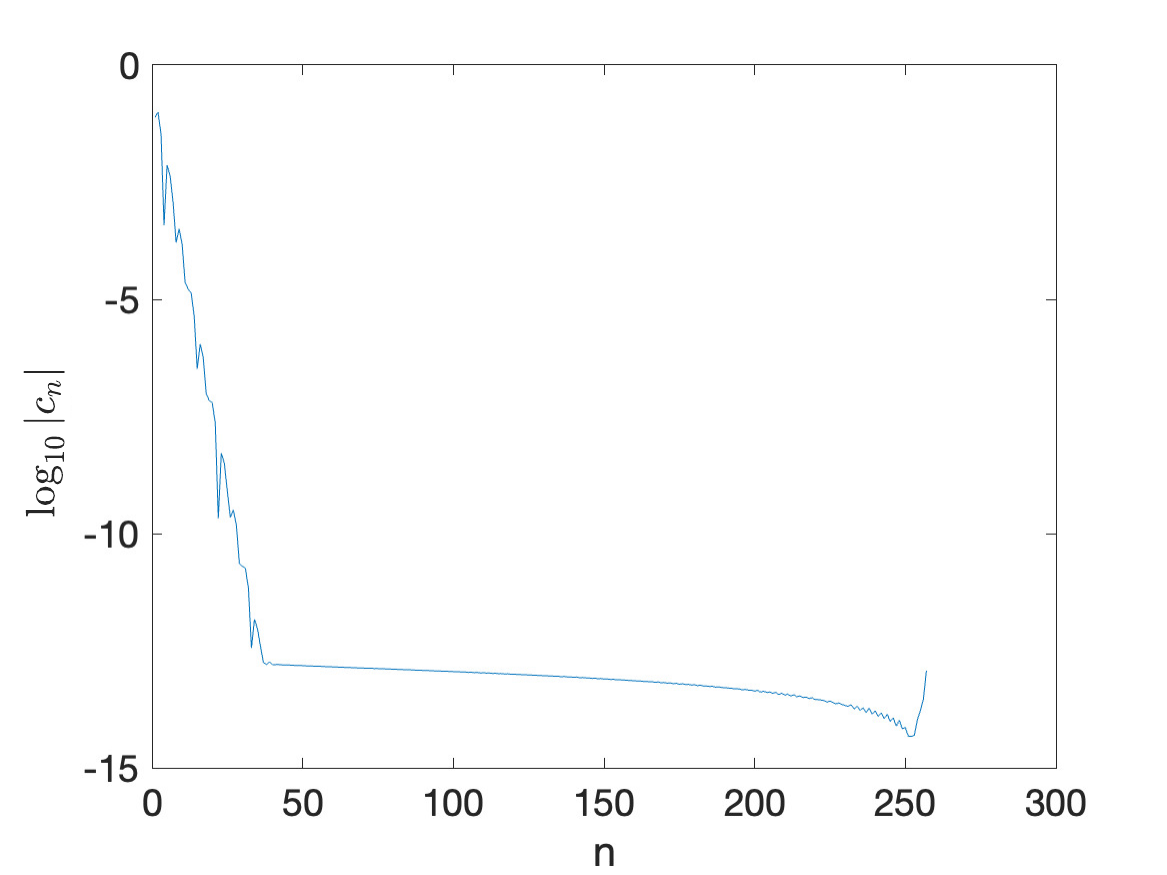}
 \includegraphics[width=0.49\textwidth]{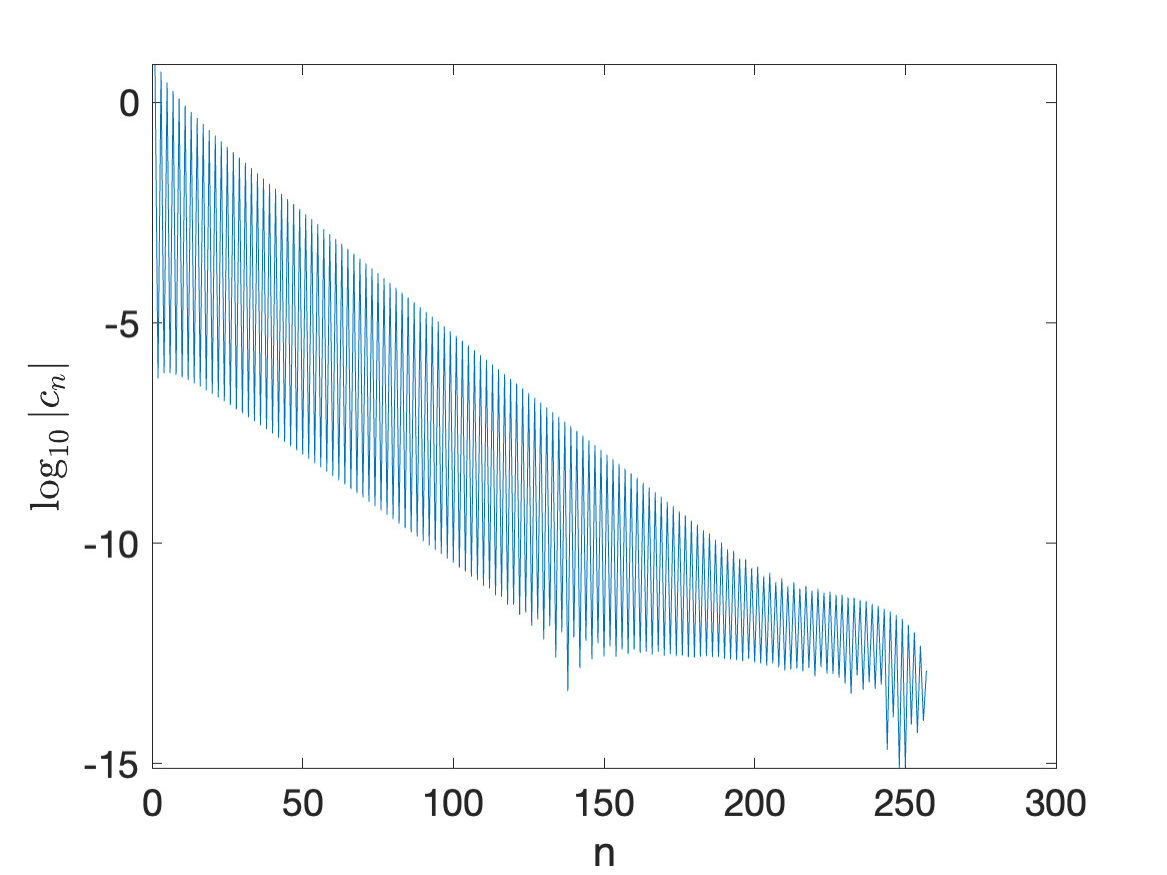}
 \caption{Chebyshev coefficients for the solitary wave of the fKdV equation for $c=1$ and 
 $\alpha=0.4$, on the left in domain I, on 
 the right in domain II.}
 \label{fKdVsol04cheb}
\end{figure}

 The Newton iteration for the multi-domain approach uses $N=256$ and 
 $b=-a=10^{-2}$ and $\delta=5*10^{-3}$. The solution is shown on the 
 right of Fig.~\ref{fKdVsol04FFT} in domain II. The Chebyshev 
 coefficients for the solution are given in Fig.~\ref{fKdVsol04cheb} 
 and indicate that the solution is well resolved. 

\begin{figure}[htb!]
 \includegraphics[width=0.49\textwidth]{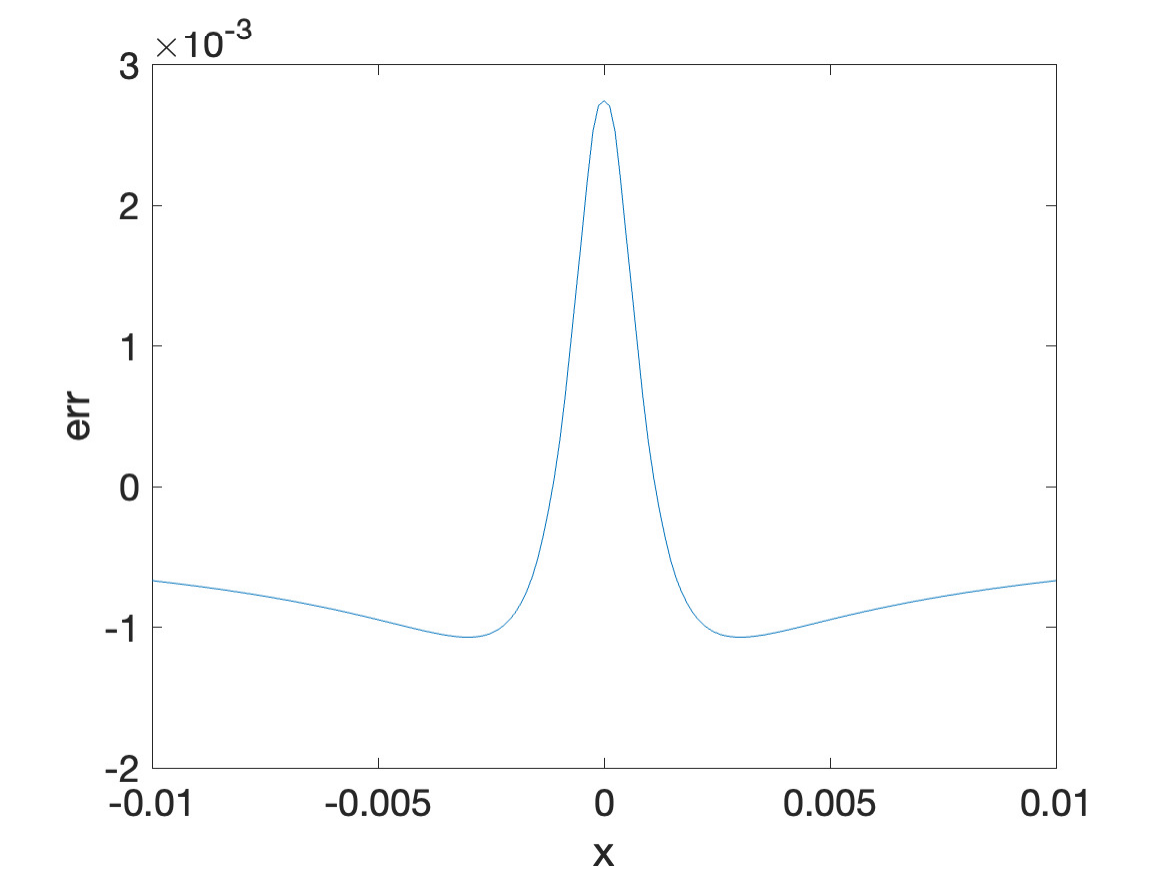}
 \includegraphics[width=0.49\textwidth]{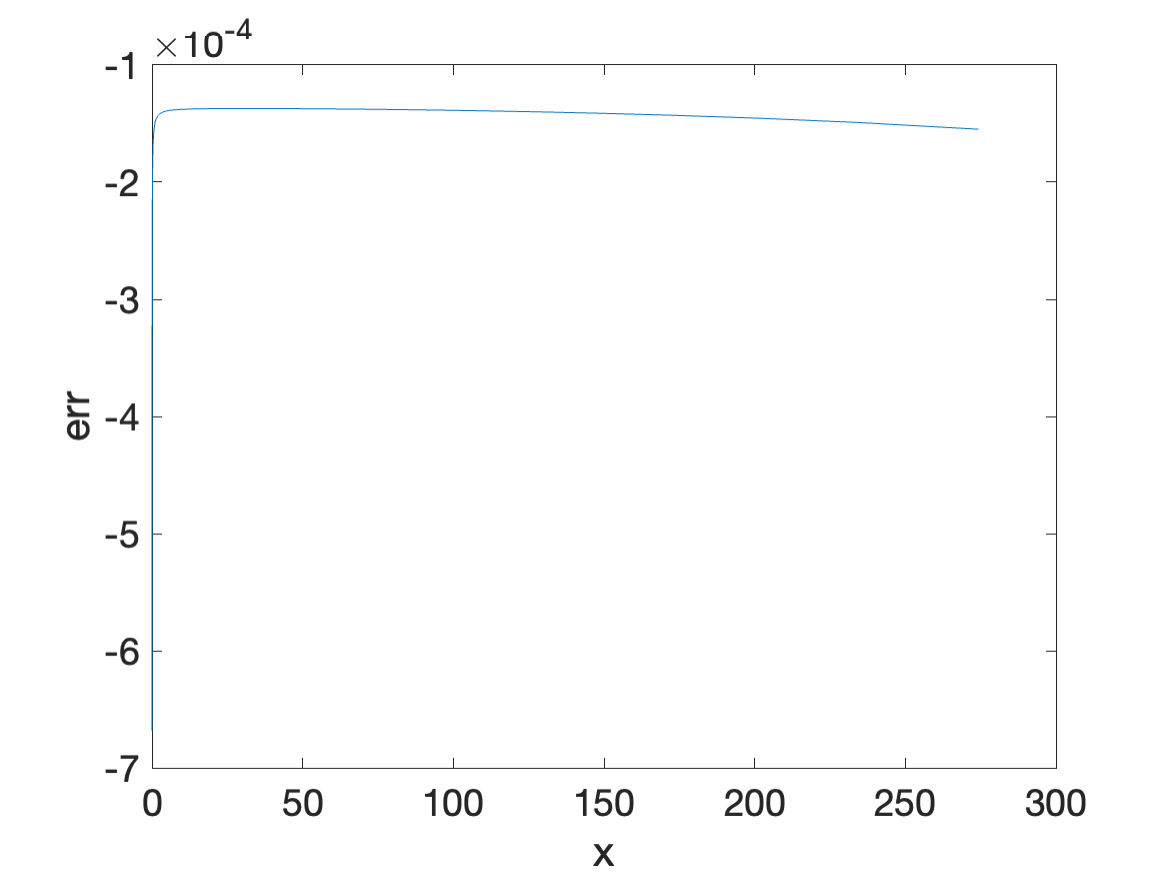}
 \caption{Difference between the solitary wave of the fKdV equation for $c=1$ and 
 $\alpha=0.4$ constructed with an FFT and a multi-domain approach for 
 various values of $x$.}
 \label{errfKdVsol04FFT}
\end{figure}

The difference between the solitary wave constructed with a 
multi-domain and with an FFT approach is shown in 
Fig.~\ref{errfKdVsol04FFT}. It is of the order of $10^{-3}$ and thus 
close to plotting accuracy. This means that the FFT approach becomes 
less and less useful for values of $\alpha$ close to the energy-critical value $1/3$. 

\begin{figure}[htb!]
 \includegraphics[width=0.49\textwidth]{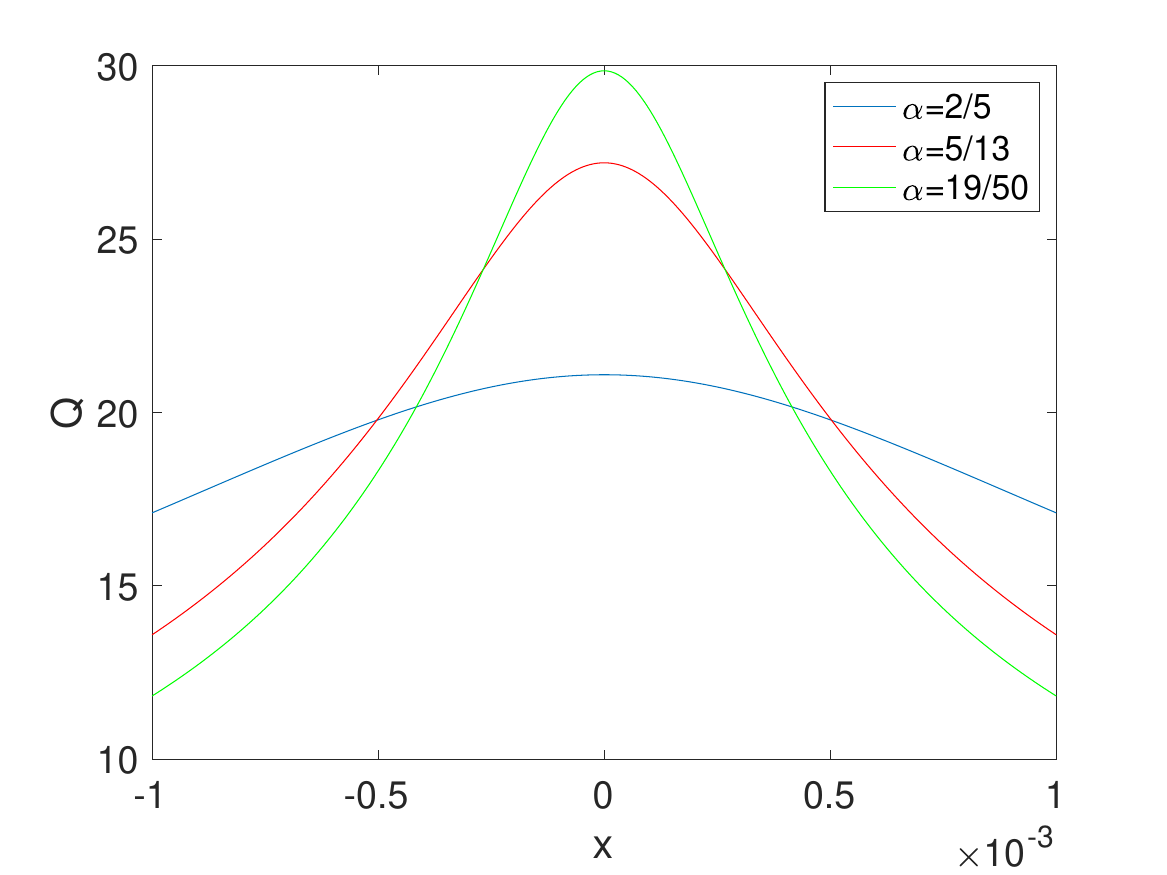}
 \includegraphics[width=0.49\textwidth]{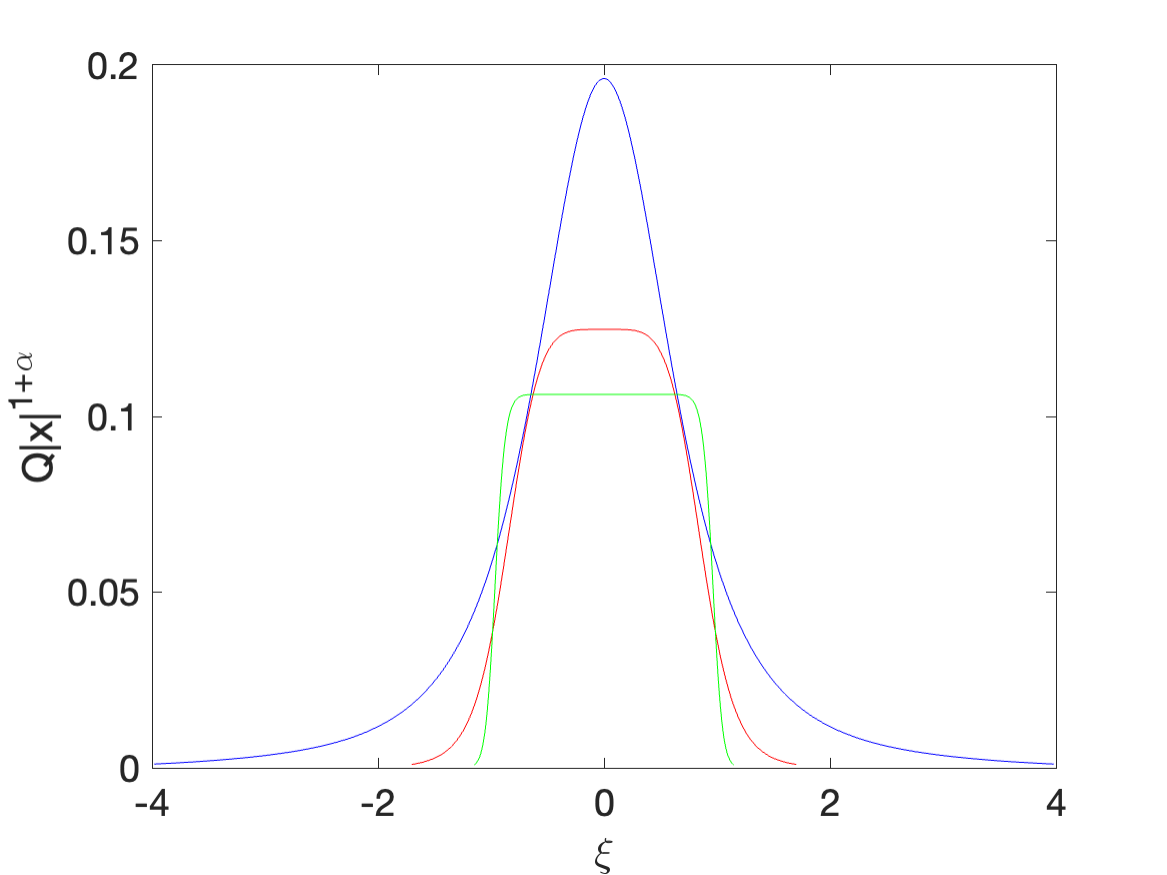}
 \caption{Solitary wave of the fKdV equation for $c=1$ and 
 $\alpha=0.4$, $\alpha=5/13$ and $\alpha=19/50$: on the left the 
 solution in domain II, on the right near infinity.}
 \label{fKdVsol513}
\end{figure}

The multi-domain approach allows to reach lower values of $\alpha$. 
For $\alpha=5/13\sim 0.385$ and $\alpha=19/50$, we use $b=-a=10^{-3}$, 
$\delta=5*10^{-4}$ and $N=400$. The solution is shown in 
Fig.~\ref{fKdVsol513} together with the solution for $\alpha=0.4$. It 
can be seen that the solution near $x=0$ has a more pronounced  peak 
closer to the origin the smaller $\alpha$ gets. The situation near 
infinity is shown on the right of the same figure. Here the solution 
becomes flatter for smaller values of $\alpha$ in the form of a 
smoothed out 
finite step. Note that we show near infinity $Q|x|^{1+\alpha}$, thus 
$Q$ seems to have some $\delta$-type behavior in the limit 
$\alpha\to1/3$. Obviously it is  numerically  challenging to 
approximate such extreme functions and to reach even 
smaller values of $\alpha$.  

\begin{figure}[htb!]
 \includegraphics[width=0.49\textwidth]{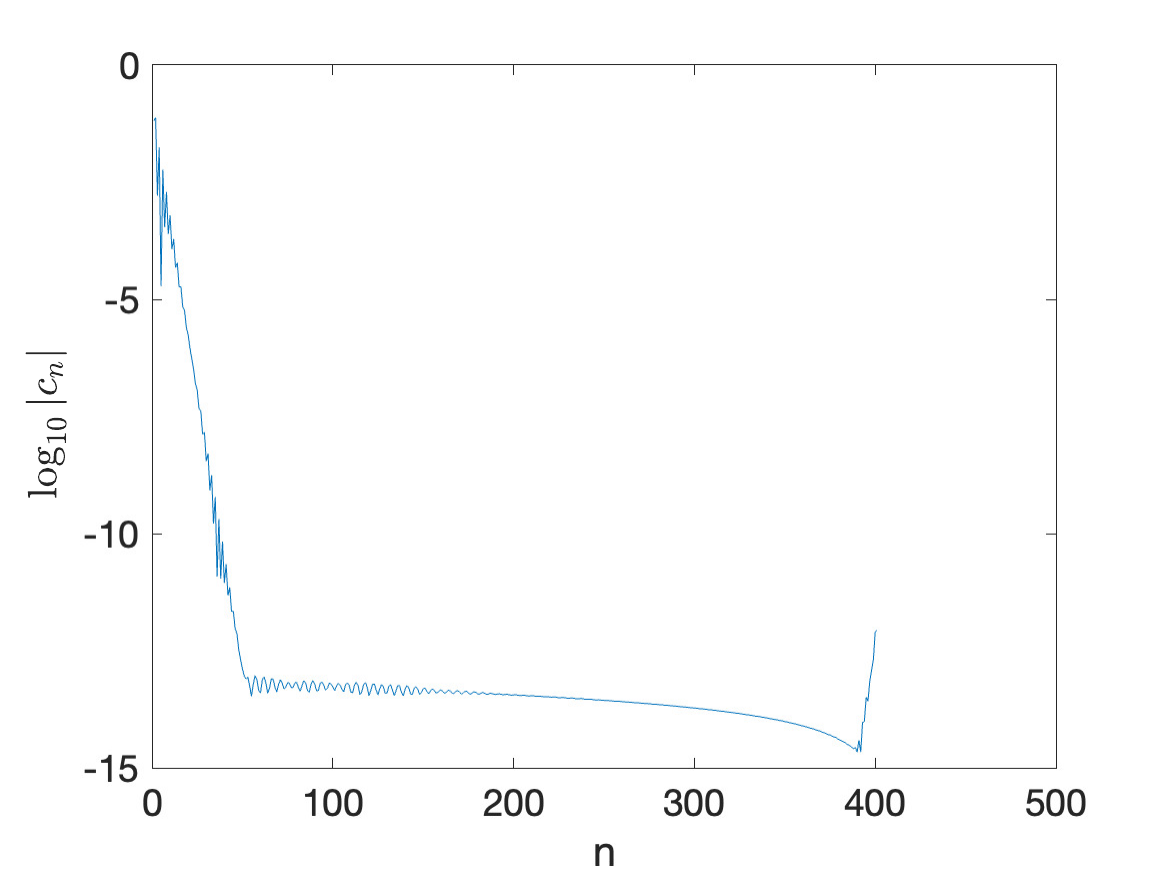}
 \includegraphics[width=0.49\textwidth]{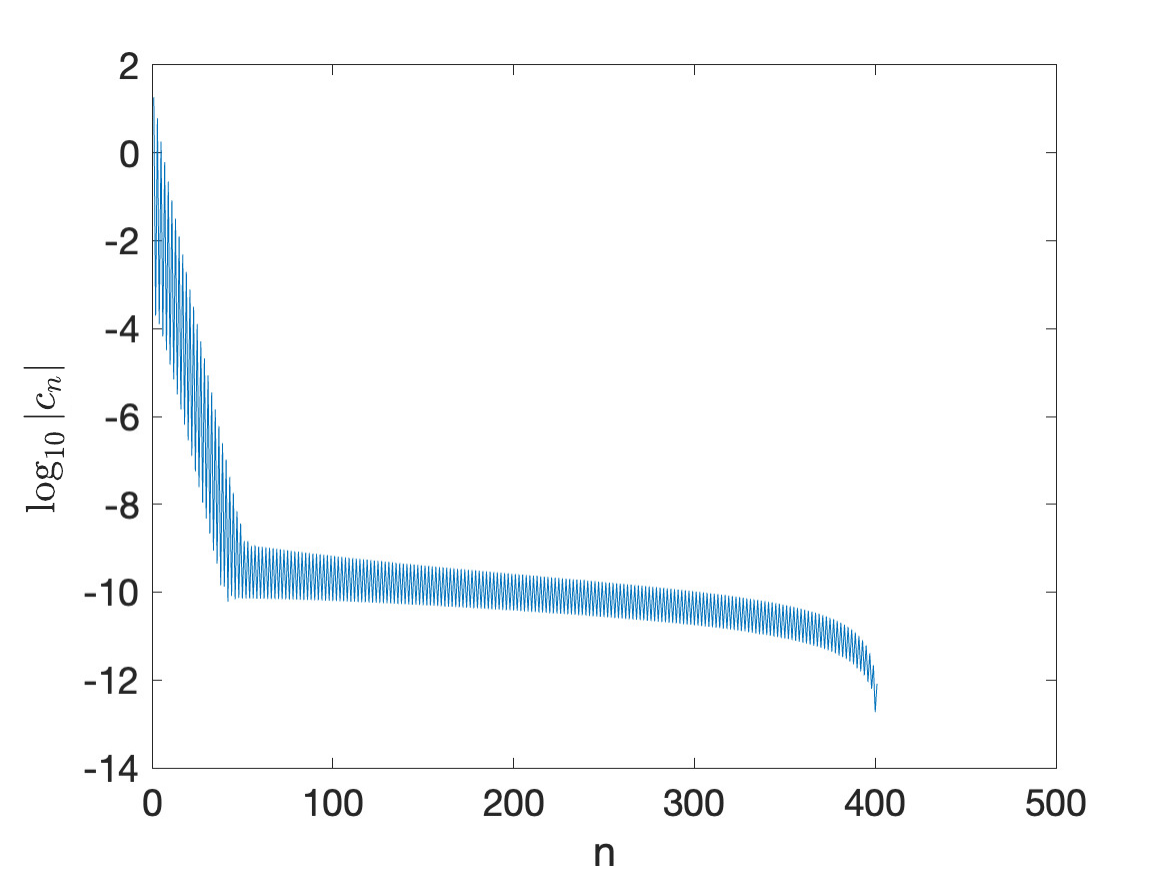}
 \caption{Chebyshev coefficients for the solitary wave of the fKdV equation for $c=1$ and 
 $\alpha=5/13$, on the left in domain I, on 
 the right in domain II.}
 \label{fKdVsol513cheb}
\end{figure}

As can be seen in Fig.~\ref{fKdVsol513cheb}, the solitary wave is 
well resolved in coefficient space for instance for $\alpha=5/13$. 
If one were interested in 
reaching even smaller values of $\alpha$, it would be necessary to 
to take even smaller finite domains, and ideally several. This is 
beyond the scope of the current paper. 

\section{Conclusion}
In this paper, we have shown that a multi-domain spectral approach 
can be efficiently used to compute the fractional Laplacian for 
rational $\alpha$ (for values of  $q$ of the order 100) to 
machine precision for smooth functions in the local parameters on the 
underlying Riemann surface. This was applied to compute solitary 
waves for the fractional KdV equation with unprecedented accuracy and 
for values closer to the critical  value than before. Note that the 
equation (\ref{Q}) also gives the solitary waves for the fractional 
modified KdV equation, \cite{KSW},
\begin{equation}
	u_{t}\pm u^{2}u_{x} - D^{\alpha}u_{x} = 0
	\label{fmKdV}
\end{equation}
and for stationary solutions $u=Q(x)e^{-ict}$ of fractional NLS 
equation \cite{KSM},
\begin{equation}
	iu_{t}-D^{\alpha}u \pm |u|^{2}u = 0.
	\label{fNLS}
\end{equation}

The defining equation for traveling wave and stationary solutions is 
after some rescaling of the form
\begin{equation}
	D^{\alpha}Q + Q -\kappa Q^{n}=0,
	\label{Q},
\end{equation}
where $\kappa>0$ is a constant (in the focusing cases, where solitons are to be 
expected, the constant is positive), and where $n\in 
\mathbb{N}$. It will be the subject of future work to explore how 
close one can get to the energy critical case in the above examples. 
Similar techniques can be applied to fractional Camassa-Holm 
equations, see \cite{KO} and references therein. In this case 
non-smooth solitons called peakons are to be expected which can be 
conveniently treated with a multi-domain approach. 

Since the previously studied FFT approaches found solutions with a 
numerical precision close to 
plotting accuracy for values of $\alpha$ close to the energy critical 
value $\alpha_{c}$, it will be interesting to study how the 
multi-domain approach will affect the time evolution of perturbed solitary waves, on the whole 
real line which has been so far mainly studied on the torus. This is 
in particular interesting for the mass critical and supercritical 
case where a blow-up in finite time is expected. This will be the 
subject of future work.

\end{document}